\documentclass[10pt]{article}

\usepackage{a4wide}
\usepackage{amssymb}
\usepackage{amsfonts}
\usepackage{amsmath}
\input xy
\xyoption{arrow} \xyoption{matrix}

\date{}

\newtheorem{proposition}{Proposition}[section]
\newtheorem{theorem}[proposition]{Theorem}
\newtheorem{lemma}[proposition]{Lemma}

\newtheorem{corollary}[proposition]{Corollary}

\def\der{\partial }

\def\nFM0{{\nu }_{F,M_0}}
\def\nFN0{{\nu }_{F,N_0}}
\def\nGN0{{\nu }_{G,N_0}}

\def\N0{ {\bf N}_0 }

\def\t{\otimes}

\def\ra{\rightarrow}

\def\Xpm{X^{\pm }}

\def\s{\sigma}
\def\Z{\mathbb{Z}}

\def\l1{{\lambda}_1}

\def\a{\alpha}
\def\a0{ {\alpha }_0}
\def\a1{ {\alpha }_1}

\def\l{\lambda}

\def\nFGM0{{\nu }_{F,G,M_0}}

\def\nFN0{{\nu}_{F,N_0}}

\def\sm{{\sigma}^m}

\def\sm1{{\sigma}^{-1}}

\def\smtp1{{\sigma}^{-t+1}}

\def\S1{S^{-1}}

\def\Xpm1{X^{\pm 1}_1}

\def\sPM1{{\sigma }^{\pm 1}}
\def\sMP1{{\sigma }^{\mp 1 }}

\def\di{{\rm d.ind}}

\def\L{\Lambda}

\def\G{\Gamma}

\def\Ytm1{Y^{t-1}}
\def\Yim1{Y^{i-1}}

\def\CL{{\cal L}}
\def\CM{{\cal M}}
\def\CN{{\cal N}}

\def\ass{{\rm ass}}
\def\bT{\overline{T}}

\def\ker{ {\rm ker } }

\def\SL2Z{ {\rm SL}_2({\bf Z}) }

\def\CR{ {\cal R}}

\def\CL{{\cal L}}

\def\Gp1{ G^{1 , 1 } }
\def\P11{ P^{-1 , 1 } }
\def\Pp1{ P^{1 , 1 } }

\def\nCLsr{{}^\nu\kern-2pt {\cal L}^{\sigma , \rho  }}
\def\nP{{}^\nu \kern-2pt P}
\def\nL{{}^\nu\kern-2pt L}
\def\nLL{{}^\nu\kern-2pt \Lambda}
\def\nPsr{{}^\nu\kern-2pt P^{\sigma , \rho  }}
\def\nLsr{{}^\nu\kern-2pt L^{\sigma , \rho  }}
\def\nuCL{{}^\nu\kern-2pt  {\cal L}}
\def\nCLsr{{}^\nu\kern-2pt {\cal L}^{\sigma , \rho  }}
\def\nCL1m{{}^\nu\kern-2pt {\cal L}^{-1 , 1  }}
\def\x1nu{x^\frac{1}{\nu}}
\def\xm1nu{x^{-\frac{1}{\nu}}}

\def\rad{{\rm rad}}

\def\CR{ {\cal R}}

\def\CN{{\cal N}}
\def\ra{\rightarrow }

\def\CB{{\cal B}}

\def\CC{ {\cal C}}

\def\nAM0{{\nu }_{{\cal A},M_0}}
\def\nAN0{{\nu }_{{\cal A},N_0}}

\def\End{ {\rm End }}

\def\CR{ {\cal R }}

\def\bR{\overline{R}}

\def\ga{\mathfrak{a}}
\def\gb{\mathfrak{b}}
\def\gc{\mathfrak{c}}

\def\gp{\mathfrak{p}}

\def\SL{{\rm SL}}

\def\Spec{{\rm Spec}}

\def\tpi{\widetilde{\pi}}

\def\di!{\frac{\der^i}{i!}}
\def\dik!{\frac{\der^k_i}{k!}}

\def\N{\mathbb{N}}

\def\0{\overline{0}}
\def\1{\overline{1}}

\def\Ln1{\L_{n,\overline{1}}}

\def\a1{a_{\overline{1}}}

\def\bs{\overline{s}}

\def\S{\Sigma}

\def\vn1{\overrightarrow{n-1}}

\def\Q{\mathbb{Q}}
\def\im{{\rm im}}

\def\mL{\mathbb{L}}

\def\gf{\mathfrak{f}}

\def\mJ{\mathbb{J}}
\def\mI{\mathbb{I}}

\def\ann{{\rm ann}}
\def\lann{{\rm l.ann}}

\def\bM{\overline{M}}

\def\K1{{\rm K}_1}

\def\hmI1{\widehat{\mI_1}}
\def\tmI1{\widetilde{\mI_1}}
\def\tmJ1{\widetilde{\mJ_1}}
\def\hB1{\widehat{B_1}}
\def\hCB1{\widehat{\CB_1}}

\def\bS{\overline{S}}

\def\Den{{\rm Den}}

\def\Ore{{\rm Ore}}

\def\Den{{\rm Den}}

\def\maxDen{{\rm max.Den}}

\def\br{\overline{r}}

\def\gt{\mathfrak{t}}
\def\bs{\overline{s}}

\def\ga{\mathfrak{a}}

\def\tor{{\rm tor}}
\def\udim{{\rm udim}}

\def\RSm{ R\langle S^{-1}\rangle}
\def\mLR{ \mathbb{L}(R)}
\def\mLlR{ \mathbb{L}_l(R)}
\def\mLrR{ \mathbb{L}_r(R)}
\def\mLsR{ \mathbb{L}_*(R)}

\def\pmLlR{ {}'\mathbb{L}_l(R)}
\def\mLrpR{ \mathbb{L}_r'(R)}

\def\mLRa{ \mathbb{L}(R, \ga )}
\def\mLlRa{ \mathbb{L}_l(R, \ga )}
\def\mLrRa{ \mathbb{L}_r(R, \ga )}
\def\mLrpRa{ \mathbb{L}_r'(R, \ga )}
\def\mLsRa{ \mathbb{L}_*(R, \ga )}

\def\pmLlR{ {}'\mathbb{L}_l(R)}
\def\mLrpR{ \mathbb{L}_r'(R)}

\def\LOresR{\mathbb{L}{\rm Ore}_*(R)}
\def\LOresRa{\mathbb{L}{\rm Ore}_*(R, \ga )}
\def\LOreR{\mathbb{L}{\rm Ore}(R)}
\def\LOrelR{\mathbb{L}{\rm Ore}_l(R)}
\def\LOrerR{\mathbb{L}{\rm Ore}_r(R)}

\def\LAOsR{\mathbb{L}{\rm AO}_*(R)}
\def\LAOlR{\mathbb{L}{\rm AO}_l(R)}
\def\LAOrR{\mathbb{L}{\rm AO}_r(R)}
\def\LAOR{\mathbb{L}{\rm AO}(R)}

\def\mL{\mathbb{L}}

\def\AOresR{{\rm AOre}_*(R)}
\def\AOrelR{{\rm AOre}_l(R)}
\def\AOrerR{{\rm AOre}_r(R)}
\def\AOreR{{\rm AOre}(R)}

\def\pCCR{{}'\CC_R}

\begin{document}

\author{V. V. \  Bavula 
}

\title{Localizable sets and the localization of a ring at a localizable set }

\maketitle

\begin{abstract}

The concepts of localizable set, localization of a ring and a module at a localizable set are introduced  and studied. Localizable sets are generalization of Ore sets and denominator sets, and the localization of a ring/module at a localizable set is a generalization of localization of a ring/module at a denominator set.  For a semiprime left Goldie ring, it is proven that the set of  maximal left localizable sets that contain all regular elements is equal to the set of maximal  left denominator sets (and they are explicitly described). For a semiprime Goldie ring, it is  proven that the following five sets coincide: the maximal   Ore sets,  the maximal  denominator sets,  the maximal left or right or two-sided  localizable sets that contain all regular elements (and they are explicitly described).  

$\noindent $



$${\bf Contents}$$
\begin{enumerate}
\item Introduction.
\item Localizable sets and the localization of a ring at a localizable set.
\item Localization of a ring at an Ore set.
 \item Localization of a ring at an almost Ore set.
\item Classification of maximal localizable sets and maximal Ore sets in semiprime Goldie ring.
 \item Localization of a module at a localizable set. 
\end{enumerate}
\end{abstract}

\section{Introduction} 

In the paper all rings are unital. When we say ring we mean a $K$-algebra over a commutative ring $K$ that belongs to the centre of the ring. 

The goal of the paper is to start to develop the most general theory of {\em one-sided fractions}. For that the following new concepts are introduced: the  almost Ore set, the  localizable set and the  localizable perfect set. Their relations are given by the chain of inclusions:
$$\{ {\rm Denominator \; sets}  \} \subseteq \{ {\rm Ore \; sets}  \} \subseteq \{ {\rm almost\;  Ore\;  sets}  \} \subseteq \{ {\rm perfect\;  localizable\;  sets }  \} \subseteq \{ {\rm  localizable\;  sets}  \} . $$

{\bf The ring $\RSm$.} Let $R$ be a ring and $S$ be a multiplicative set in $R$ (that is $SS\subseteq S$, $1\in S$ and $0\not\in S$). Let $R\langle X_S\rangle$ be a ring freely generated by the ring $R$ and a set $X_S=\{ x_s\, | \, s\in S\}$ of free noncommutative indeterminates (indexed by the elements of the set $S$). Let us consider the factor ring

\begin{equation}\label{RbSbm}
\RSm := R\langle X_S\rangle/ I_S
\end{equation}
of the ring $R\langle X_S\rangle$ at the ideal $I_S$ generated by the set of elements $\{ sx_s-1, x_ss-1 \, | \, s\in S\}$. 

The kernel of the ring homomorphism
\begin{equation}\label{RbSbm1}
R\ra \RSm , \;\; r\mapsto r+ I_S
\end{equation}
is denoted by $\ass (S) = \ass_R(S)$. The  ideal $\ass_R(S)$ of $R$ has a complex structure, its description is given in Proposition \ref{B19Jan19} when $S$ is a left localizable set.  Lemma \ref{b15Jan19},  Proposition \ref{A12Jan19}.(1) and  its proof  describe a large chunk of the ideal $\ass_R(S)$,  which is the ideal $\ga (S)$. The proof of Proposition \ref{A12Jan19} contains an explicit description of the ideal $\ga (S)$.  The ideal $\ga (S)$ is the key part in the definition of perfect localizable sets. \\

{\bf Localizable sets.} 

{\it Definition.} A multiplicative set $S$ of a ring $R$ is called a {\em left localizable set} of $R$ if  
$$\RSm = \{ \bs^{-1} \br \, | \, \bs \in \bS, \br \in \bR\}\neq \{ 0\}$$
where $\bR = R/ \ga$,  $\ga = \ass_R(S)$ and $\bS = (S+\ga ) / \ga$, i.e., every element of the ring $\RSm$ is a left fraction $\bs^{-1} \br$ for some elements  $\bs \in \bS$ and $ \br \in \bR$. Similarly,  a multiplicative set $S$ of a ring $R$ is called a {\em right  localizable set}  of $R$  if 
$$\RSm = \{  \br \bs^{-1}\, | \, \bs \in \bS, \br \in \bR\}\neq \{ 0\}, $$
 i.e., every element of the ring $\RSm$ is a right  fraction $ \br\bs^{-1}$ for some elements  $\bs \in \bS$ and $ \br \in \bR$. A right and left localizable set of $R$ is called a {\em localizable set} of $R$.\\
 
  The sets of left localizable, right localizable and localizable sets of $R$ are denoted by $\mLlR$, $\mLrR$ and $\mLR$, respectively. Clearly, $ \mLR = \mLlR \cap \mLrR$.  In order to work with these three sets simultaneously we use the following notation $\mLsR$  where $*\in \{ l,r, \emptyset \}$  and $\emptyset$ 
is the empty set   $(\mL (R) = \mL_\emptyset (R))$.
 Let 
\begin{equation}\label{xLzRa}
 \ass \, \mL_* (R) =\{ \ass_R (S) \, | \, S\in \mL_*(R) \}.
\end{equation} 
 For an ideal $\ga $ of $R$, let $ \mLsRa =\{ S\in \mLsR \, | \, \ass_R(S) = \ga \}$. Then 

\begin{equation}\label{LzRa}
\mLsR= \coprod_{\ga \in \ass\,  \mLsR } \mLsRa
\end{equation}
is a disjoint union of non-empty sets. \\

{\bf The ideals $\ga (S)$, ${}'\ga (S)$ and $\ga'(S)$.} For each element $r\in R$, let $r\cdot : R\ra R$, $x\mapsto rx$ and $\cdot r : R\ra R$, $x\mapsto xr$. The sets $\pCCR := \{ r\in R\, | \, \ker (\cdot r)=0\}$ and $\CC_R' := \{ r\in R\, | \, \ker (r\cdot )=0\}$ are called the {\em sets of left and right regular elements} of $R$, respectively.  Their intersection $\CC_R=\pCCR \cap \CC_R'$ is the {\em set of regular elements} of $R$. The rings $Q_{l,cl}(R):= \CC_R^{-1}R$ and $Q_{r,cl}(R):= R\CC_R^{-1}$ are called the {\em classical left and right quotient rings} of $R$, respectively. Goldie's Theorem states that the ring $Q_{l, cl}(R)$ is  a semisimple Artinian ring iff the ring $R$ is  semiprime, $\udim (R)<\infty$ and the ring $R$ satisfies the a.c.c. on left annihilators ($\udim$ stands for the uniform dimension).

\begin{proposition}\label{A12Jan19}
Let $R$ be a ring and $S$ be a non-empty subset of $R$.
\begin{enumerate}
\item Suppose that there exists an ideal $\gb$ of $R$  such that $(S+\gb ) / \gb \subseteq \CC_{R/ \gb}$. Then there is  the least ideal, say $\ga= \ga (S)$, that satisfies this property. 
\item Suppose that there exists an ideal $\gb$ of $R$  such that $(S+\gb ) / \gb \subseteq {}'\CC_{R/ \gb}$. Then there is  the least ideal, say ${}'\ga= {}'\ga (S)$, that satisfies this property;  and $ {}'\ga (S)\subseteq \ga (S)$.
\item Suppose that there exists an ideal $\gb$ of $R$  such that $(S+\gb ) / \gb \subseteq \CC_{R/ \gb}'$. Then there is  the least ideal, say $\ga'= \ga' (S)$, that satisfies this property;  and $\ga' (S)\subseteq \ga (S)$.
\end{enumerate}
\end{proposition}

For a multiplicative set $S$ in a ring $R$, we fix the following notation (unless it is stated otherwise):  ${}'\ga = {}'\ga (S)$ and $\ga' = \ga' (S)$ (see Proposition \ref{A12Jan19}), 
\begin{equation}\label{pRpS1}
{}'R:=R/ {}'\ga{}, \;\; \pi : R\ra {}'R,\;\;   r\mapsto {}'r= r+{}'\ga , \;\; {}'S= {}'\pi (S) ,
\end{equation}
\begin{equation}\label{pRpS2}
R':= R/ \ga'; \;\; \pi' : R\ra R', \;\;  r\mapsto r'= r+\ga', \;\; S'=\pi' (S).
\end{equation}

The proof of Proposition \ref{A12Jan19} is given in Section \ref{LOCSET}. The ideals $\ga (S)$, ${}'\ga (S)$ and $\ga' (S)$ are defined in an explicit way, see (\ref{ala}), (\ref{ala1}) and (\ref{ala2}), respectively. They play an important role in the proofs of many results of this paper.

\begin{lemma}\label{b15Jan19}
Given $S\in \mLsR$ where $*\in \{ l,r, \emptyset \}$. Then $\ass_R (S) \supseteq \ga (S)$ where $\ga (S)$ is the least ideal of $R$ such that $(S+\ga (S) ) / \ga (S) \subseteq \CC_{R/\ga (S)}$, see Proposition \ref{A12Jan19}.(1).
\end{lemma}

{\bf The structure of the ring $\RSm$ and its universal property.}  Let $R$ be a ring. A 
multiplicative subset $S$ of $R$   is called  a {\em
left Ore set} if it satisfies the {\em left Ore condition}: for
each $r\in R$ and
 $s\in S$, $ Sr\bigcap Rs\neq \emptyset $.
Let $\Ore_l(R)$ be the set of all left Ore sets of $R$.
  For  $S\in \Ore_l(R)$, $\ass_l (S) :=\{ r\in
R\, | \, sr=0 \;\; {\rm for\;  some}\;\; s\in S\}$  is an ideal of
the ring $R$. 


A left Ore set $S$ is called a {\em left denominator set} of the
ring $R$ if $rs=0$ for some elements $ r\in R$ and $s\in S$ implies
$tr=0$ for some element $t\in S$, i.e., $r\in \ass_l (S)$. Let
$\Den_l(R)$ be the set of all left denominator sets of $R$. For
$S\in \Den_l(R)$, let $S^{-1}R=\{ s^{-1}r\, | \, s\in S, r\in R\}$
be the {\em left localization} of the ring $R$ at $S$ (the {\em
left quotient ring} of $R$ at $S$). Let us stress that in Ore's method of localization one can localize {\em precisely} at left denominator sets.
 In a similar way, right Ore and right denominator sets are defined. 
Let $\Ore_r(R)$ and 
$\Den_r(R)$ be the set of all right  Ore and  right   denominator sets of $R$, respectively.  For $S\in \Ore_r(R)$, the set  $\ass_r(S):=\{ r\in R\, | \, rs=0$ for some $s\in S\}$ is an ideal of $R$. For
$S\in \Den_r(R)$,  $RS^{-1}=\{ rs^{-1}\, | \, s\in S, r\in R\}$   is the {\em right localization} of the ring $R$ at $S$. 

Given ring homomorphisms $\nu_A: R\ra A$ and $\nu_B :R\ra B$. A ring homomorphism $f:A\ra B$ is called an $R$-{\em homomorphism} if $\nu_B= f\nu_A$.  A left and right  set is called an {\em Ore set}.  Let $\Ore (R)$ and 
$\Den (R)$ be the set of all   Ore and    denominator sets of $R$, respectively. For
$S\in \Den (R)$, $$S^{-1}R\simeq RS^{-1}$$ (an $R$-isomorphism)
 is  the {\em  localization} of the ring $R$ at $S$, and $\ass (R):=\ass_l(R) = \ass_r(R)$.

For a ring $R$ and   $*\in \{ l,r, \emptyset \}$, $\Den_*(R, 0)$ be the set of  $*$ denominator sets $T$ of $R$ such that $T\subseteq \CC_R$, i.e., the multiplicative set $T$ is a $*$ Ore set of $R$ that consists of regular elements of the ring $R$. For a ring $R$, we denote by $R^\times$ its  group of units (invertible elements) of the ring $R$. Theorem \ref{L11Jan19} describes the structure and the universal property of the ring $\RSm$.  

\begin{theorem}\label{L11Jan19}
Let $S\in \mLsRa$ where  $*\in \{ l, \emptyset \}$, $\bR = R/ \ga$, $\pi : R\ra \bR$, $r\mapsto \br = r+\ga$ and $\bS = \pi (S)$. Then 
\begin{enumerate}
\item $\bS\in \Den_*(\bR , 0)$.
\item The ring $\RSm$ is $R$-isomorphic to the ring $\bS^{-1}\bR$. 
\item Let $\gb$ be an ideal of $R$ and $\pi^\dag : R\ra R^\dag =R/\gb$, $r\mapsto r^\dag =r+\gb$. If $S^\dag =\pi^\dag  (S)\in \Den_*(R^\dag , 0)$ then $\ga \subseteq \gb$ and the map $\bS^{-1}\bR \ra {S^\dag}^{-1}R^\dag $, $\bs^{-1}\br \mapsto {s^\dag}^{-1}r^\dag $ is a ring epimorphism with kernel $\bS^{-1}(\gb/\ga)$. So, the ideal $\ga$ is the least ideal $\ga$ of the ring $R$ such that $S+\ga \in \Den_*(R/\ga  , 0)$.
\item Let $f: R\ra Q$ be a ring homomorphism such that $f(S)\subseteq Q^\times$ and the ring $Q$ is generated by $f(R)$ and the set $\{ f(s)^{-1}\, | \, s\in S\}$. Then 
\begin{enumerate}
\item $\ga \subseteq \ker (f)$ and the map $\bS^{-1}\bR \ra Q$, $\bs^{-1}\br\mapsto f(s)^{-1} f(r)$ is a ring epimorphism with kernel $\bS^{-1} (\ker (f) / \ga )$, and $Q= \{ f(s)^{-1} f(r)\, | \, s\in S, r\in R\}$. 
\item Let $\widetilde{R}=R/\ker (f)$ and $\widetilde{\pi}:R\ra \widetilde{R}$, $r\mapsto \tilde{r}= r+\ker (f)$. Then $\widetilde{S}:=\widetilde{\pi}(S)\in \Den_l(\widetilde{R}, 0)$ and $\widetilde{S}^{-1}\widetilde{R}\simeq Q$, an $\widetilde{R}$-isomorphism. 
\end{enumerate}
\end{enumerate}
\end{theorem}  

A similar result holds for $*=r$, i.e., for right localizable sets. Statements 3 and 4 of Theorem \ref{L11Jan19} are the universal property of localization of a ring at a (left or right) localizable set. In the particular case when $S\in \Den_*(R)$, these are precisely the universal property of localization at a (left or right) denominator set. \\

 In view of Theorem \ref{L11Jan19}.(1,2), for $S\in \mLsR$ we denote by $S^{-1}R$ the ring $\RSm$  for $*\in \{ l,\emptyset\}$ and by $RS^{-1}$ for $*\in \{ r,\emptyset\}$. In particular, for $S\in \mLR$, 
 $\RSm = S^{-1}R\simeq RS^{-1}$. Elements of the rings  $S^{-1}R$ and $ RS^{-1}$ are denoted by $s^{-1}r$ and $rs^{-1}$, respectively, where $s\in S$ and $r\in R$. \\

{\bf Perfect localizable sets.} By Lemma \ref{b15Jan19}, $\ass_R(S) \supseteq \ga (S)$ for all $S\in \mLsR$ where $*\in \{ l,r, \emptyset \}$.\\ 

{\it Definition.} A localizable set $S\in \mLsR$ is called {\em perfect} if $\ass_R(S) = \ga (S)$, i.e., the ideal $\ass_R(S)$ is `the least possible'.\\

 Therefore, localizations at perfect localizable sets are the most free/largest possible. Another  feature of perfect localizable sets is that the ideal $\ass_R(S) = \ga (S)$ admits an explicit description that can be computed  in many examples (see the proof of Proposition \ref{A12Jan19} and (\ref{ala})). \\

Let $\mL_*^p(R)=\{ S\in \mLsR\, | \, \ass_R(S) = \ga (S)\}$ and $\ass\, \mL_*^p(R)=\{ \ass_R(S) \,  | \, S\in \mL_*^p(R)\}$. Clearly, 
\begin{equation}\label{LpRU}
\mL_*^p(R)=\coprod_{\ga \in \ass\, \mL_*^p(R)}\mL_*^p(R, \ga )
\end{equation}
where $\mL_*^p(R, \ga ) =\{ S\in \mL_*^p(R)\, | \, \ass_R(S) = \ga \}$. Clearly, $\mL_*^p(R, \ga ) = \mL_*^p(R)\cap \mL_*(R, \ga )$. \\

{\bf The sets ${}'\mLlR$,  $\mL_r'(R)$ and ${}'\mL_{l,r}'(R)$.} We denote by $\Ore_*(R)$ (where $*\in \{ l,r, \emptyset \}$) the set of $*$ Ore sets of  $R$. So, $\Ore_l(R)$ is the set of all left Ore sets of $R$. \\

{\it Definition.} Let ${}'\mLlR$ (resp.,   $\mL_r'(R)$) be the set of all multiplicative sets $S$ of $R$ such that ${}'S\in \Ore_l({}'R )$ (resp., $ S'\in \Ore_r(R')$) (see (\ref{pRpS1}) and (\ref{pRpS2})). Let  
 $${}'\mL_{l,r}'(R):={}'\mLlR\cap \mL_r'(R)\;\; {\rm and}\;\;  \LOresR:=\mLsR \cap \Ore_*(R)$$ where $*\in \{ l,r, \emptyset \}$. The elements of the set $ \LOresR$ are called the {\em localizable  $*$  Ore sets}. \\

{\bf Localizable left/right Ore sets.}   The study of localizations at left, right, and left and right  Ore sets was started in the paper \cite{Bav-intdifline}. In particular, 
\cite[Theorem 4.15]{Bav-intdifline} states that {\em every Ore set is a localizable set}, i.e.,   $\Ore (R)\subseteq \mLR$. Therefore, 
\begin{equation}\label{LOre=Ore}
\LOreR=\Ore (R).
\end{equation}
This fact also follows from Theorem \ref{10Jan19} (see Theorem \ref{16Jan19} for details). Proposition \ref{C15Jan19} establishes relations between the concepts that are introduced above.

\begin{proposition}\label{C15Jan19}

\begin{enumerate}
\item $\LOrelR \subseteq \pmLlR \subseteq \mL^p_l(R)$ and  $\LOrelR = \pmLlR \cap \Ore_l(R)= \mL^p_l(R)\cap \Ore_l(R)$.
\item $\LOrerR \subseteq \mLrpR \subseteq \mL^p_r(R)$ and  $\LOrerR = \mLrpR \cap \Ore_r(R)= \mL^p_r(R)\cap \Ore_r(R)$.
\item  $\Ore (R) \subseteq {}'\mL_{l,r}'(R) \subseteq \mL^p (R)$. 
\end{enumerate}
\end{proposition}

{\bf Criterion for a left Ore set to be a left localizable set.} Theorem \ref {B16Jan19} is such a criterion.  

\begin{theorem}\label{B16Jan19}
Let $R$ be a ring, $S\in \Ore_l(R)$, and  $\ga = \ass_R(S)$. Then
\begin{enumerate}
\item $S\in \mLlR$ iff ${}'\ga \neq R$ where the ideal ${}'\ga = {}'\ga (S)$ of $R$ is as in Proposition \ref{A12Jan19}.(2) and (\ref{ala1}). 
\item  Suppose that ${}'\ga \neq R$. Let ${}'\pi : R\ra {}'
R:=R/ {}'\ga$, $r\mapsto {}'r= r+{}'\ga$ and ${}'S={}'\pi (S)$. Then 
\begin{enumerate}
\item $'S\in {}'\Den_l({}'R)$.
\item $\ga = {}'\pi^{-1}(\ass_l({}'S))$.
\item $S^{-1}R\simeq {}'S^{-1}{}'R$, an $R$-isomorphism. 
\end{enumerate}
\end{enumerate}
\end{theorem}
Theorem \ref{C17Jan19} is a criterion for a right Ore set to be a localizable set. \\

{\bf Localization at an Ore set.} 	Theorem \ref{10Jan19} is the reason why every Ore set is localizable.

\begin{theorem}\label{10Jan19}
Let $R$ be a ring and $S\in \Ore (R)$.
\begin{enumerate}
\item $\ga :=\{ r\in R\, | \, srt=0$ for some elements $s,t\in S\}$ is an ideal of $R$ such that $\ga \neq R$. 
\item Let $\pi : R\ra \bR :=R/\ga$, $r\mapsto \br =r+\ga$. Then $\bS :=
\pi (S) \in \Den (\bR , 0)$, $\ga = \ga (S)=\ass_R(S)$ and $S^{-1}R\simeq \bS^{-1}\bR$, an $R$-isomorphism. 
\item Let $\gb$ be an ideal of $R$ and $\pi^\dag :R\ra R^\dag :=R/\gb$, $r\mapsto r^\dag =r+\gb$. If $S^\dag := \pi^\dag (S)\in \Den (R^\dag ,0)$ then $\ga \subseteq \ga^\dag $ and the map $\bS^{-1}\bR \ra {S^\dag}^{-1}R^\dag $, $\bs^{-1}\br \mapsto {s^\dag}^{-1}r^\dag$ is a ring epimorphism. 
\item Let $f: R\ra Q$ be a ring homomorphism such that $f(S)\subseteq Q^\times$ and the ring $Q$ is generated by $f(R)$ and $\{ f(s)^{-1}\, | \, s\in S\}$. Then  
\begin{enumerate}
\item $\ga \subseteq \ker (f)$ and the map $\bS^{-1}\bR \ra Q$, $ \bs^{-1}\br \mapsto f(s)^{-1}f(r)$ is a ring epimorphism with kernel $\bS^{-1}(\ker (f) / \ga )$. 
\item Let $\widetilde{R} = R/\ker (f)$ and $\tpi : R\ra \widetilde{R}$, $r\mapsto r+\ker (f)$. Then $\widetilde{S}:= \tpi (S) \in  \Den (\widetilde{R}, 0)$ and $\widetilde{S}^{-1}\widetilde{R}\simeq Q$, an $\widetilde{R}$-isomorphism. 
\end{enumerate}
\end{enumerate} 
\end{theorem} 

Theorem \ref{10Jan19}.(1,2) states that for every  $S\in \Ore (R)$, the ideal $\ga$ in Theorem \ref{10Jan19}.(1) coincides with the ideals $\ga (S)=\ass_R(S)$. So, every Ore set $S$ of $R$ are localizable and  and the localization $S^{-1}R$ of the ring $R$ at the localizable set $S$ is $R$-isomorphic to the localization $\bS^{-1}\bR$ of the ring $\bR$ at the denominator set $\bS$ of $\bR$ (see Theorem \ref{10Jan19} and Theorem \ref{16Jan19}).

For a ring $R$ and and its ideal $\ga$, let 
\begin{eqnarray*}
 {}'\Den_l(R , \ga)&:=&\{ S\in \Den_l(R)\, | \, \ass_l(S) = \ga, S\subseteq \pCCR \}, \\
\Den_r'(R, \ga )&:=&\{ S\in \Den_r(R)\, | \, \ass_r(S) = \ga, S\subseteq \CC_R' \}.
\end{eqnarray*}

Corollary \ref{c10Jan19} offers  a view from another angle at the localization of a ring  at an Ore set. 

\begin{corollary}\label{c10Jan19}
Let $R$ be a ring, $S\in \Ore (R)$ and $\ga = \ass_R(S)$. We keep the notation of Theorem \ref{10Jan19}. Let $\ga_l:= \ass_l(S)$ and $\pi_l: R\ra R_l:=R/\ga_l$, $r\mapsto r+\ga_l$; $\ga_r:= \ass_r(S)$ and $ \pi_r:R\ra R_r:= R/\ga_r$, $r\mapsto r+\ga_r$. Then  
\begin{enumerate}
\item $\ga_l+\ga_r\subseteq \ga$.
\item $S_l:= \pi_l(S)\in \Den_r'(R_l, \ga / \ga_l)$ and $R_lS_l^{-1}\simeq S^{-1}R\simeq RS^{-1}$, $R$-isomorphisms. 
\item $S_r:= \pi_r(S)\in  {}'\Den_l(R_r, \ga / \ga_r)$   and $S_r^{-1}R_r\simeq S^{-1}R\simeq RS^{-1}\simeq R_lS_l^{-1}$, $R$-isomorphisms. 
\item ${}'\ga (S) = \ga_r$ and $ \ga'(S) = \ga_l$. 
\end{enumerate}
\end{corollary}

Proposition \ref{A13Jan19} explains the origin of the construction of localizable sets. For every Ore set it gives an explicit construction of a {\em localizable set}  with the same localization.  For a ring $R$ and and its ideal $\ga$, let $ \Ore (R, \ga ):=\{ S\in \Ore (R)\, | \, \ass_R (S)=\ga \}$.

\begin{proposition}\label{A13Jan19}
Let $R$ be a ring, $S\in \Ore (R, \ga )$, $\pi : R\ra \bR := R/ \ga$, $R\mapsto \br = r+\ga$ and $\widetilde{S}= S+\ga$. Then $S\subseteq \widetilde{S}$, $\widetilde{S}\in \mLRa$ and $\widetilde{S}^{-1}R\simeq S^{-1}R$, an $R$-isomorphism. 
\end{proposition}
 
{\bf The set $\max  \mLsR $ maximal elements in $\mLsR$ where $*\in \{ l,r,\emptyset\}$.}  For a ring $R$, the set $\maxDen_l(R)$ of maximal left denominator sets  (w.r.t. $\subseteq$) is a {\em non-empty}
 set, \cite[Lemma 3.7.(2)]{larglquot}. Let $\max \, \mLsR$ be the set of maximal elements (w.r.t. $\subseteq $) of the set $\mLsR$.

  \begin{theorem}\label{C12Jan19}
Let $R$ be a ring. Then $\max  \mLsR\neq \emptyset$.
\end{theorem}
 
 The key idea of the proof of Theorem \ref{C12Jan19} is to use Lemma \ref{a12Jan19} and Zorn's Lemma. 
 
 \begin{lemma}\label{a12Jan19}
Let $R$ be a ring, $S\in \mLsRa$ and $T\in \mL_*(R, \gb )$ such that $S\subseteq T$ where $*\in \{ l,r,\emptyset\}$. Then $\ga \subseteq \gb$ and for  $*\in \{ l,\emptyset\}$  the map $S^{-1}R\ra T^{-1}R$, $s^{-1}r\mapsto t^{-1}r$ is an $R$-homomorphism with kernel $S^{-1} (\gb / \ga ) = \bS^{-1}(\gb / \ga )$ where  $\bS = \{ s+\ga \, | \, s\in S\}$. A similar result holds for $*=r$. 
\end{lemma}

  {\bf Classification of maximal Ore sets of a semiprime  Goldie ring.}  It was proved that the set $\maxDen_l(R)$ is a finite set if the classical left quotient ring $Q_{ l, cl}(R):=\CC_R^{-1}R$ of $R$ is
 a semisimple  Artinian ring, \cite{Bav-Crit-S-Simp-lQuot},  or a left Artinian ring,  \cite{Bav-LocArtRing}, or a left Noetherian ring, \cite{Crit-lNoeth-lQuot}. In each of the three cases an explicit description of the set $\maxDen_l(R)$ is given.  For a ring $R$, let $\min (R)$ be the set of its minimal prime ideals. For a ring $R$, the rings $Q_{l, cl}(R):=\CC_R^{-1}R$ and $Q_{r, cl}(R):=R\CC_R^{-1}$ are called the {\em classical left and right quotient rings}, respectively. 
 
 The next theorem is an explicit description of maximal Ore sets of a semiprime Goldie ring.
 
\begin{theorem}\label{17Jan19}
Let $R$ be a semiprime Goldie ring and $\CN_* :=\{ S\in \max  \mL_*(R) \, | \, \CC_R\subseteq S\}$ where $*\in \{ l,r,\emptyset\}$. Then 
\begin{enumerate}
\item $\max  \Ore (R)= \max \Den (R) = \{ \CC (\gp ) \, | \, \gp \in\min (R)\}=\CN_*$ for all $*\in \{ l,r,\emptyset\}$  where $\CC (\gp ):= \{ c\in R \, | \, c+\gp \in \CC_{R/\gp} \}$. So, every maximal Ore set of $R$ is a maximal denominator set, and vice versa. 
\item For all $S\in \max  \Ore (R)$ the ring $S^{-1}R$ is a simple Artinian ring. 
\item $Q_{cl} (R) \simeq \prod_{S\in \max  \Ore (R)} S^{-1}R$ (where $Q_{cl} (R) := \CC_R^{-1}R\simeq  R\CC_R^{-1}$ is the classical quotient ring of $R$).
\item $\max  \ass\, \Ore (R) = \ass \, \max  \Ore (R)=\min (R)$. In particular, the ideals in the set $\ass \, \max \, \Ore (R)$ are incomparable. 
\end{enumerate}
\end{theorem}

{\bf Classification of maximal left localizable sets of a semiprime left Goldie ring that contain the set of regular elements of the ring.} Theorem \ref{18Jan19} is such a classification. 

\begin{theorem}\label{18Jan19}
Let $R$ be a semiprime left Goldie ring and $\CN :=\{ S\in \max  \mLlR \, | \, \CC_R\subseteq S\}$. Then $\CN = \max  \Den_l(R)=\{ \CC (\gp ) \, | \, \gp \in \min (R)\}$. 
\end{theorem}

So, every  maximal left localizable set of a semiprime left Goldie ring that contain the set of regular elements of the ring is a maximal left denominator set, and vice versa. The proof of Theorem \ref{18Jan19} is based on Theorem \ref{A18Jan19}. 

\begin{theorem}\label{A18Jan19}
Let $R$ be a ring and $S_1, \ldots , S_n\in \mLsR$ where $*\in \{ l,r, \emptyset \}$, $\gp_i=\ass_R(S_i)$, $R_i:= R/\gp_i$ and $Q_i$ be the localization of $R$ at $S_i$. Suppose that the rings $Q_i$ are simple Artinian rings, $\bigcap_{i=1}^n\gp_i=0$ and  $\bigcap_{j\neq i}^n\gp_j\neq 0$ for $i=1, \ldots , n$. Then 
\begin{enumerate}
\item The rings $R_i$ are semiprime $*$ Goldie rings.
\item $\min (R) = \{ \gp_1, \ldots , \gp_n\}$. 
\item $Q_i\simeq Q_{*, cl}(R_i)$ for $i=1, \ldots , n$ (an $R$-isomorphism).
\item $S_i\subseteq \CC (\gp_i )$ for $i=1, \ldots , n$. 
\item $Q_{*, cl}(R)\simeq \prod_{i=1}^n Q_i$. 
\item For all $i=1, \ldots , n$, $\CC (\gp_i) \in \max \, \mLsR$. 
\end{enumerate}
\end{theorem}




\section{Localizable sets and the localization of a ring at a localizable set}\label{LOCSET}

In this section, (left, right) localizable sets of a ring and the construction of localization of a ring at them are introduced and studied.  A proof of Theorem \ref{L11Jan19} is given. A criterion for a left (resp., right) Ore set to be a left (resp., right)  localizable set is presented,  Theorem \ref{B16Jan19}.(1) (resp., Theorem \ref{C17Jan19}.(1)).  A description of the ideal $\ass_R(S)$ of $R$ is given, Proposition \ref{B19Jan19}. Proofs of Proposition \ref{A12Jan19},  Proposition \ref{C15Jan19}, Theorem \ref{C12Jan19}, Lemma \ref{a12Jan19},  Theorem \ref{18Jan19} and Theorem \ref{A18Jan19} are given. 

The following notation is fixed: $R$ is a ring, $S$ is a multiplicative  set of $R$, $\ass_l(S):=\{ r\in R\, | \, sr=0$ for some $s\in S\}$ and $\ass_r(S):=\{ r\in R\, | \, rs=0$ for some $s\in S\}$.   We use standard terminology on localizations of a ring at denominator sets, see \cite{Jategaonkar-LocNRings,MR,Stenstrom-RingQuot}.\\

{\bf Proof of Proposition  \ref{A12Jan19}.} We  keep the notation of  Proposition   \ref{A12Jan19}. Let $\G$ be the set of ordinals. The ideal  $\ga$ (resp., ${}'\ga$, $\ga'$) is the union 
\begin{equation}\label{aappa}
\ga= \bigcup_{\l \in \G}\ga_\l\;\;\;   ({\rm resp.,}\;\;  {}'\ga = \bigcup_{\l \in \G}{}'\ga_\l, \;\; \ga' = \bigcup_{\l \in \G}\ga_\l') 
\end{equation}
of ascending chain of ideals $\{ \ga_\l \}_{\l \in \G}$ (resp., $\{ {}'\ga_\l  \}_{\l \in \G}$, $\{ \ga'_\l \}_{\l \in \G}$), where $\l \leq \mu$ in $\G$ implies $\ga_\l \subseteq \ga_\mu$ (resp., ${}'\ga_\l \subseteq {}'\ga_\mu$, $\ga_\l' \subseteq \ga_\mu'$). The ideal $\ga_\l$ (resp., ${}'\ga_\l$, $\ga_\l'$) are defined inductively as follows: the ideals $\ga_0=\ga (S, R)$ (resp.,  ${}'\ga_0$,  $\ga_0'$) is generated by the set $\{ r\in R\,  | \, sr=0$ or $ rt=0$ for some elements $s,t\in S\}$ (resp., $\{ r\in R\,  | \,  rt=0$ for some element $t\in S\}$, $\{ r\in R\,  | \, sr=0$ for some element $s\in S\}$), and  for $\l \in \G$ such that $\l >0$ (where below $\Big( \{ \ldots \} \Big)$ means the ideal of $R$ generated by the set  $ \{ \ldots \} $), 
\begin{equation}\label{ala}
\ga_\l = \begin{cases}
 \bigcup_{\mu <\l \in \G}\ga_\mu & \text{if $\l$ is a limit ordinal},\\
\Big( \{ r\in R\, | \, sr\in \ga_{\l -1}\; {\rm or}\; rt\in \ga_{\l -1} \; {\rm for \; some} \; s,t\in S\}\Big) & \text{if if $\l$ is not a limit ordinal}.\\
\end{cases}
\end{equation}
(resp., 
\begin{equation}\label{ala1}
{}'\ga_\l = \begin{cases}
 \bigcup_{\mu <\l \in \G}{}'\ga_\mu & \text{if $\l$ is a limit ordinal},\\
\Big( \{ r\in R\, | \, rt\in {}'\ga_{\l -1}\; {\rm for \; some} \; t\in S\}\Big) & \text{if if $\l$ is not a limit ordinal}.\\
\end{cases}
\end{equation}

\begin{equation}\label{ala2}
\ga_\l' = \begin{cases}
 \bigcup_{\mu <\l \in \G}\ga_\mu' & \text{if $\l$ is a limit ordinal},\\
\Big( \{ r\in R\, | \, sr\in \ga_{\l -1}'\; {\rm for \; some} \; s\in S\}\Big)  & \text{if if $\l$ is not a limit ordinal}.)\;\;\; \Box\\
\end{cases}
\end{equation}

The first ordinal such that the ascending chain of ideals in (\ref{aappa}) stabilizes is called the $\G$-{\em length} of the ideal denoted $l_\G (\ga )$ (resp., $l_\G ({}'\ga )$, $l_\G (\ga' )$).
 
\begin{lemma}\label{a15Jan19}
Let $R$ be a ring, $S$ be a multiplicative set, $\ga = \ass_R(S)$, $\gb$ be an ideal of the ring $R$ such that $\gb \subseteq \ga$, $\widetilde{R}=R/\gb$ and $\widetilde{S}= (S+\gb ) / \gb= \{ s\in \gb \, | \, s\in S\}$. Then 

\begin{enumerate}
\item $\widetilde{S}$ is a multiplicative set of the ring $\widetilde{R}$, $\widetilde{R}\langle \widetilde{S}^{-1}\rangle \simeq \RSm$ and $\ass_{\widetilde{R}}(\widetilde{S})=\ga / \gb$. 
\item Let $\bR = R/ \ga$ and $\bS = (S+\ga ) / \ga$. Then $\RSm \simeq \bR \langle \bS^{-1}\rangle$, $\bS\subseteq \CC_{\bR}$ and $\ass_{\bR}(\bS ) = \ga / \ga =0$. 
\end{enumerate}
\end{lemma}

{\it Proof}. Straightforward. $\Box $\\

{\bf Proof of Lemma \ref{b15Jan19}.} Let $\ga = \ass_R(S)$ and $\bR = R/ \ga$.  Then $\bS = (S+\ga ) / \ga \in \CC_{\bR}$. Hence, $\ga \supseteq \ga (S)$, by the minimality of the ideal $\ga (S)$, see Proposition \ref{A12Jan19}.(1). $\Box$\\

For a ring $R$, its ideal $\ga$,  and   $*\in \{ l,r, \emptyset \}$, let $\Den_*(R, \ga )$ be the set of  $*$ denominator sets $S$ of $R$ such that $\ass_*(S)=\ga $. 

Lemma \ref{c15Jan19}  and Lemma \ref{d15Jan19} show that denominator sets are localizable sets, the localization of a ring at a  denominator set is the same as the localization of a ring at the denominator set treated as a localizable set, and $\ass_*(S) = \ga (S)$ for all $S\in \Den_*(R, \ga )$. So, the localization at a localizable set is a generalization of the localization an a denominator set.

\begin{lemma}\label{c15Jan19}
Given $S\in \Den_*(R, \ga )$ where $*\in \{ l,r,\emptyset\}$. Then $\ga = \ga (S)$ where $\ga (S)$ is the least ideal of $R$ such that $(S+\ga (S) ) / \ga (S)\in \CC_{R/ \ga (S)}$, see Proposition \ref{A12Jan19}.(1). 
\end{lemma}

{\it Proof}. By the very definition of the ideal $\ga$, $\ga \subseteq \ga (S)$. Since $(S+\ga ) / \ga \in \CC_{R/\ga }$, we have the inverse inclusion $\ga \supseteq \ga (S)$, by the minimality of the ideal $\ga (S)$. Therefore  $\ga = \ga (S)$. 
 $\Box $

\begin{lemma}\label{d15Jan19}
$\Den_*(R)\subseteq \mLsR$ where $*\in \{ l,r,\emptyset\}$. For all ideals $\ga $ of $R$, $\Den_*(R, \ga )\subseteq \mLsRa$ and for all $S\in \Den_l(R, \ga )$, $S^{-1}R\simeq R\langle S^{-1}\rangle$ and $\ga = \ga (S)=\ass_R(S)$. 
\end{lemma}

{\it Proof}. Given $S\in \Den_l(R, \ga )$. By Lemma \ref{b15Jan19} and Lemma \ref{c15Jan19}, $\ga = \ga (S) \subseteq \ass_R(R)$. Let $\bR = R/ \ga $ and $\bS = (S+\ga ) / \ga$. Then 
$$\bS \in \Den_l(R, 0)\;\; {\rm  and}\;\; S^{-1}R\simeq \bS^{-1}\bR .$$ The ring $\bS^{-1} \bR$ is generated by the ring $\bR$ and the set $\{ \bs^{-1} \, | \, \bs \in \bS \}$. Therefore, there is natural ring epimorphism $\RSm\ra \bS^{-1}\bR$ which is an $R$-homomorphism. Hence, $\ass_R(S) \subseteq \ga$, and so $\ga = \ga (S)=\ass_R(S)$ and $\RSm \simeq \bS^{-1}\bR$.  $\Box $\\

{\bf Proof of Theorem \ref{L11Jan19}.} To prove the theorem it suffices to consider the case $*=l$. 

1 and 2. Let $\CR = \RSm$. By Lemma \ref{a15Jan19}.(2), 
$$\RSm \simeq \bR \langle \bS^{-1}\rangle, \;\; \bS\subseteq \CC_{\bR}\;\; {\rm  and}\;\; \ass_{\bR}(\bS )=0.$$
 Hence, $\bR \subseteq \CR$. Then, for all elements $\bs \in \bS$ and $\br \in \bR$, 
 $$\br \bs^{-1} = \bs_1^{-1}\br_1^{-1}$$ for some elements $\bs_1\in \bS$ and $\br_1\in \bR$. Then $\bs_1\br = \br_1\bs $ in $\bR$. Hence, $\bS \in \Den_l(\bR , 0)$, and so $\bR \langle \bS^{-1}\rangle \simeq \bS^{-1} \bR$, by Lemma \ref{d15Jan19}.

3. There is a natural $R$-epimorphism $\RSm\ra {S^\dag}^{-1}R^\dag$, and statement 3 follows from statement 2. 

4. The ring homomorphism $r: R\ra Q$ induces the ring epimorphism $f: \RSm\ra Q$ since the ring $Q$ is generated by $f(R)$ and the set $\{ f(s)^{-1}\, | \, s\in S\}$. Since $\RSm \simeq \bS^{-1} \bR$ (statement 2), the statement (a) follows.  

In particular, $Q=\{ f(s)^{-1} f(r) \, | \, s\in S, r\in R\}$. Furthermore, $Q=\{ f(\bs)^{-1} f(\br) \, | \, s\in S, r\in R\}$  since $\ga \subseteq  \ker (f)$. Hence,  $Q=\{ f(\tilde{s})^{-1} f(\tilde{r}) \, | \,  s\in S, r\in R\}$. Therefore, $\widetilde{S} \in \Den_l(\widetilde{R}, 0)$ and $\widetilde{S}^{-1}\widetilde{R}\simeq Q$, via $f$. $\Box$

\begin{corollary}\label{a16Jan19}
If $S\in \mLRa$ then $S\in \mLlRa$, $S\in \mLrRa$ and the localizations of $R$ at $S$ as a localizable set, a left localizable set and a right localizable set are $R$-isomorphic. 
\end{corollary}

{\it Proof}. The corollary follows from Theorem \ref{L11Jan19}.(1,2).  $\Box $\\

{\bf Proof of Lemma \ref{a12Jan19}.} Recall that  $\ga = \ass_R(S)$ and $\gb = \ass_R(T)$. Let $Q$ be a  subring of $T^{-1}R$ which is generated by the images of the ring $R$ and the set $\{ s^{-1} \, | \, s\in S\}$ in $T^{-1}R$ (recall that $S\subseteq T$). Applying Theorem \ref{L11Jan19}.(4a) 
to the ring homomorphism $R\ra Q\subseteq T^{-1}R$, $ r\mapsto \frac{r}{1}$ we obtain the ring $R$-homomorphism 
$$S^{-1}R\ra T^{-1}R, \;\; s^{-1}r\mapsto s^{-1} r.$$
 Since 
 $S^{-1}R= \bS^{-1}\bR$ and $T^{-1}R= \bT^{-1}(R/ \gb )$ where $\bT = \{ t+\gb \, | \, t\in T\}$, the kernel of the $R$-homomorphism is $ \bS^{-1} (\gb / \ga )$. $\Box$ \\

{\bf The maximal elements in $\mLsR$ where $*\in \{ l, r, \emptyset \}$.} 

\begin{lemma}\label{b12Jan19}
Let $R$ be a ring,     $*\in \{ l,r, \emptyset\}$,  and  a set $\{ S_i\}_{i\in I}\subseteq \mLsR$ be such that for any two elements  $S_i$ and $S_j$ there is an element $S_k$ such that $S_i\cup S_j\subseteq S_k$. Then $S:= \bigcup_{i\in I}S_i\in \mLsR$ and the ring $S^{-1}R =\injlim S_i^{-1}R  $ is the injective limit of $R$-homomorphism of rings  $\{ S_i\}_{i\in I}$ given by the $R$-homomorphisms $S_i^{-1}R\ra S_j^{-1}R$ in Lemma \ref{a12Jan19} in case $S_i\subseteq S_j$. 
\end{lemma}

{\it Proof}. Straightforward (use Lemma \ref{a12Jan19}).  $\Box $\\

{\bf Proof of Theorem \ref{C12Jan19}.} The Theorem \ref{C12Jan19}  follows from Zorn's Lemma and Lemma \ref{b12Jan19}. $\Box$ \\

{\bf Classification of maximal left localizable sets of a semiprime left Goldie ring that contain the set of regular elements.} \\

{\bf Proof of Theorem \ref{A18Jan19}.} 1.  By the assumption $\bigcap_{i=1}^n\gp_i =0$. So, we have   ring monomorphisms
$$ R\ra \CR := \prod_{i=1}^n R_i\ra Q:= \prod_{i=1}^n Q_i, \;\; r\mapsto (r+\gp_1, \ldots , r+\gp_n), \;\; (r_1, \ldots , r_n) \mapsto \bigg(\frac{r_1}{1}, \ldots , \frac{r_n}{1}\bigg).$$
We  identify the rings $R$ and $\CR$ with their images in the ring $Q$. 

(i) $R_i$ {\em is a prime $*$ Goldie ring with $Q_{*, cl} (R_i)\simeq Q_i$, an $R$-isomorphism}: By Theorem \ref{L11Jan19}.(1,2), $$\bS_i:=(S_i+\gp_i)/\gp_i\in \Den_*(R_i, 0)\;\; {\rm  and}\;\; \bS_i^{-1}R_i\simeq Q_i,$$
 an $R$-isomorphism. Since the ring $Q_i$ is a simple Artinian ring, $\bS_i\subseteq \CC_{R_i} \subseteq Q_i^\times$. Hence, $\CC_{R_i}\in \Den_*(R, 0)$ and $Q_i\simeq Q_{*, cl} (R_i)$. By Goldie's Theorem, $R_i$ is a prime $*$ Goldie ring. 

(ii) $\gp_i\in \Spec (R)$ {\em for} $i=1, \ldots , n$: By the statement (i), the ring $R_i$ is prime, and so  $\gp_i\in \Spec (R)$. 

(iii) {\em The ring $R$ is a semiprime $*$ Goldie ring}: The ring $R$ is a semiprime ring since $\bigcap_{i=1}^n \gp_i =0$ and the ideals $\gp_i$ are prime. 

Let $\udim_{*, R}$ denote the $*$ uniform dimension of the $*$ $R$-module $(\udim_{l,R}$ is the left uniform dimension, $\udim_{r,R}$ is the right  uniform dimension and $\udim_R$ stands for $\udim_{l,R}$ and $\udim_{r,R}$). The ring $R_i=R/\gp_i$ is a prime $*$ Goldie ring. Hence, $\udim_{*,R}(R_i) = \udim_{*,R}(R_i)<\infty$. Since $R\subseteq \prod_{i=1}^n R_i$, 
$$ \udim_{*, R} (R) \leq \udim_{*,R}\bigg(\prod_{i=1}^n R_i\bigg)=\sum_{i=1}^n \udim_{*, R}(R_i)=\sum_{i=1}^n \udim_{*, R_i}(R_i)<\infty.$$
Let $X$ be a non-empty subset of $R$ and $*.\ann_R(X)$ be its $*$ annihilator ($\lann_R(X)=\{ r\in R\, | \, rX=0\}$ is the left annihilator of $X$ in $R$, etc). Since $R\subseteq \CR$, $*.\ann_R(X)= R\cap *.\ann_\CR (X)$. By the statement (i), the ring $\CR = \prod_{i=1}^n R_i$ satisfies the a.c.c. on $*$ annihilators, hence so does the ring $R$. The proof of the statement (iii) is complete. 

2. (i) $\min (R) \subseteq \{ \gp_1, \ldots , \gp_n\}$: Given $\gp \in \min (R)$. Then $\bigcap_{i=1}^n \gp_i=\{ 0\} \subseteq \gp$,hence $\gp_i\subseteq \gp$ for some $i$, and so $\gp_i = \gp$, by the minimality of $\gp$. 

(ii)  $\min (R) = \{ \gp_1, \ldots , \gp_n\}$: By the statement (i), the ring $R$ is semiprime, i.e., $\bigcap_{\gp \in \min (R)}\gp =0$. Since $\min (R) \subseteq \{ \gp_1, \ldots , \gp_n\}$ and $\bigcap_{j\neq i} \gp_i\neq 0$ for all $i=1, \ldots , n$, we must have $\min (R) = \{ \gp_1, \ldots , \gp_n\}$.

3. Statement 3 has already been proven, see the statement (i) in the proof of statement 1. 

4. See the proof of the statement (i) in the proof of statement 1. 

5. The ring $R$ is a semiprime $*$ Goldie ring. By \cite[Theorem 4.1]{Bav-Crit-S-Simp-lQuot}, 
$$Q_{*, cl}(R)\simeq \prod_{\gp \in \min (R)}Q_{*, cl}(R/\gp ).$$ By statement 2, $\min (R)= \{ \gp_1, \ldots , \gp_n\}$ and $Q_i\simeq Q_{*, cl}(R_i)$ for $i=1, \ldots , n$ (statement 3). Now, statement 5 follows.

6. By \cite[Theorem 4.1]{Bav-Crit-S-Simp-lQuot}, 
$$\CC (\gp_i) \in \max  \Den_*(R), \;\; \ass_* \, \CC (\gp_i) = \gp_i\;\;  {\rm and }\;\;\CC(\gp_i)^{-1} R\simeq Q_i,$$ a simple Artinian ring. Suppose that $\CC (\gp_i) \subseteq T$ for some $T\in \mLsR$. By Lemma \ref{a12Jan19}, $$\gp_i=\ass_R(\CC (\gp_i))\subseteq \ga := \ass_R(T)$$ and there is an $R$-homomorphism $Q_i= \CC (\gp_i )^{-1}R\ra R\langle T^{-1}\rangle$. Since $Q_i$ is a simple Artinian ring, $\gp_i = \ga$. Then, by Theorem \ref{L11Jan19}.(1), $\bT= (T+\gp_i)/\gp_i\in \Den_*(R_i, 0)$. Hence, $\bT \subseteq \CC_{R_i}$, and so $T\subseteq \CC (\gp_i)$, and statement 6 follows. $\Box$\\

{\bf Proof of Theorem \ref{18Jan19}.} Let $\CM = \max \, \Den_l(R)$. By \cite[Theorem 4.1]{Bav-Crit-S-Simp-lQuot}, $$\CM = \{ \CC (\gp ) \, | \, \gp \in \min (R)\}$$ and the set $\CM$ satisfies the conditions of Theorem \ref{A18Jan19}.  By Theorem \ref{A18Jan19}.(2,6), $$\CM \subseteq \max\, \mLsR.$$
 Since $\CC_R\subseteq \CC (\gp )$ for all $\gp \in \min (R)$, $\CM \subseteq \CN$. To finish the proof it remains to show that $\CN \subseteq \CM$. Given $S\in \CN $. Let $\ga = \ass_R(S)$.

(i) $\ga = \bigcap_{\gp \in \min (R), \ga\subseteq \gp} \gp$:  Since $S\in \CN$, $\CC_R\subseteq S$,  and so there is an $R$-homomorphisms $f: Q:=\CC_R^{-1}R\ra S^{-1}R$. Therefore, $\ga = R\cap \ker (f)$ and $\ker (f) = \CC_R^{-1} \ga$. By \cite[Theorem 4.1]{Bav-Crit-S-Simp-lQuot}, 
$$\CC_R^{-1}R\simeq \prod_{\gp \in \min (R)} \CC(\gp )^{-1}R$$ is a finite direct product of simple Artinian rings $\CC(\gp )^{-1}R$ with $\ass_R(\CC (\gp )) = \gp$, and the statement (i) follows. 

(ii) $S\subseteq \CC (\gp )$ {\em for some (unique) $\gp \in \min (R)$ such that} $\ga \subseteq \gp$: By Theorem \ref{L11Jan19}.(1), $\bS := (S+\ga ) / \ga \in \Den_l(\bR , 0)$ where $\bR = R/\ga$ is a semiprime ring with $$\min (\bR ) = \{ \gp / \ga \, | \, \gp \in\min (R), \ga \subseteq \gp \},$$ by the statement (i). Since $\CC_R\subseteq \CC (\gp )$ for all $\gp \in \min (R)$, we have that $ \overline{\CC}_R:= (\CC_R+\ga ) / \ga \subseteq\CC_{\bR}$ since 
$$\bR = R/\ga = R/\bigcap_{\gp \in D}\gp \subseteq \prod_{\gp \in D}R/ \gp$$ where
$D:= \{ \gp \in \min (R)\, | \, \ga\subseteq \gp \}$, see the statement (i). 
 Since $\CC_R\in \Ore_l(R)$, $ \overline{\CC}_R\in \Den_l(\bR , 0)$ and $$ \overline{\CC}_R^{-1}\bR \simeq \CC_R^{-1}(R/ \ga ) \simeq \CC_R^{-1}R/\CC_R^{-1}\ga \simeq \prod_{\gp \in D}\CC(\gp )^{-1}R.$$
 Therefore, $\bR$ is a semiprime left Goldie ring with $\min (\bR ) = \{ \gp / \ga \, | \, \gp \in D\}$.
 
 Since $\bS \in \Den_l(\bR , 0)$, the left denominator set $\bS$ is contained in a maximal denominator set $\bR$, i.e., $\bS \subseteq \CC_{\bR}(\gp / \ga )$ for some prime ideal $\gp \in D$. Then $S\subseteq \CC_R(\gp )$ since $ \bR /(\gp / \ga )\simeq R/ \gp $. $\Box$ \\
 
 {\bf The set $T_{*,\ga}(R)$ and the ring $Q_{*, \ga} (R)$.} Lemma \ref{b17Jan19}.(4) describes the largest element $T_{*,\ga}(R)$ in the set $\mLsRa$.

 \begin{lemma}\label{b17Jan19}
Let $R$ be a ring and $\ga \in \ass \, \mLsR$ where $*\in \{ l,r,\emptyset\}$. 
\begin{enumerate}
\item The set $\mLsRa$ is a commutative multiplicative semigroup where for $S, T \in \mLsRa$, $ST$ is the submonoid of $(R,\cdot )$ generated by $S$ and $T$. 
\item For all $S, T \in \mLsRa$, $S\cup T\subseteq ST$.
\item For all $S_1, S_2, T\in \mLsRa$ such that $S_1\subseteq S_2$, $TS_1\subseteq TS_2$. 
\item The set $T_{*, \ga} :=\bigcup_{S\in \mLsRa}S$ is the largest element in the set $\mLsRa$ (w.r.t. $\subseteq $), and $ST_{*, \ga }= T_{*, \ga }$ for all $S\in \mLsRa$. 
\item For all $S, T \in \mLsRa$ such that  $S\subseteq T$, $S^{-1}R\subseteq  T^{-1}R\subseteq T_{*, \ga}(R)^{-1}R= \bigcup_{S'\in \mLsRa}S'^{-1}R$. 
\end{enumerate}
\end{lemma}

{\it Proof}. Let $\pi : R\ra \bR := R/ \ga$, $r\mapsto \br = r+\ga$.

1. By Theorem \ref{10Jan19}.(1), $\bS , \bT \in \Den_*(\bR , 0)$. By \cite[Theorem 2.1.(1)]{larglquot}, $\overline{ST} \in \Den_*(\bR , 0)$. Hence, $\bS \bT = \bS \, \bT \in \Den_*(\bR , 0)$ and $\overline{ST}\subseteq ((\overline{ST})^{-1}\bR)^\times$. By Theorem \ref{L11Jan19}.(3), $ST\in \mLsRa$. 

2 and 3. Statements 2 and 3 are obvious.

4. Statement 4 follows from statement 2 and Lemma \ref{b12Jan19}.

5. Statement 5 follows from  Lemma \ref{a12Jan19} and  Lemma \ref{b12Jan19}. $\Box $\\

{\it Definition.} Let $R$ be a ring and $\ga \in \ass \, \mLsRa$ where $*\in \{ l,r,\emptyset\}$. The localization of the ring $R$ at $T_{*, \ga } (R)$ (see Lemma \ref{b17Jan19}.(5) is denoted by  $\Q_{*, \ga } (R)$. \\


A left or right Ore set of a ring is called {\em regular} if it consists of regular elements of the ring. A regular left/right Ore set is automatically a left/right denominator set. 
\begin{theorem}\label{3Jul10}
\cite[Theorem 2.1]{larglquot}  For a ring $R$ there is  a largest (w.r.t.
 inclusion) regular left  Ore set $S_l(R)$ in $R$ and  the ring 
 $Q_l(R ):= S_l(R)^{-1}R$ is called the {\em largest left quotient ring} of  $R$.
\end{theorem}

 For a ring $R$ there is  a largest regular   Ore set $S(R)$ in $R$ and  the ring 
 $Q (R ):= S(R)^{-1}R\simeq RS(R)^{-1}$ is called the {\em largest  quotient ring} of  $R$, \cite[Theorem 4.1.(2)]{larglquot}. The interested reader is referred to \cite{larglquot} for more information about the largest regular left/right Ore sets and the largest left/right quotient rings. Notice that for a ring $R$, its classical left/right quotient ring does not always exists.

Theorem \ref{B17Jan19}, describes the multiplicative sets  $T_{*, \ga } (R)$ and the rings $\Q_{*, \ga } (R)$ via the  largest regular Ores sets and the largest quotient rings.

\begin{theorem}\label{B17Jan19}
Let $R$ be a ring, $\ga \in \ass \, \ass \, \mL_*(R)$ and  $\pi : R\ra \bR := R/\ga$, $r\mapsto \br = r+\ga$. Then 
\begin{enumerate}
\item $T_{*, \ga } (R)=\pi^{-1} (S_*(\bR ))$ where $S_*(\bR )$  is the largest regular $*$ Ore set in $\bR$. 
\item $\Q_{*, \ga } (R) \simeq Q_*(\bR )$, an $R$-isomorphism  where $Q_*(\bR)= S_*(\bR)^{-1}\bR$ for $*\in \{ l, \emptyset \}$ and $Q_r(\bR)= \bR S_r(\bR)^{-1}$. 
\item $T_{*, \ga } (R)=\s^{-1} (\Q_{*, \ga } (R)^\times )=\s^{-1} (Q_*(\bR )^\times )$ where $\s : R\ra \Q_{*, \ga } (R) \simeq Q_*(\bR )$, $r\mapsto \frac{r}{1}$.
\end{enumerate}
\end{theorem}

{\it Proof}. 1. The set $T:= \pi^{-1} (S_*(\bR ))$ is a multiplicative set in $R$ since $\pi (T) = S_*(\bR )$ is so. Clearly, $T_*:= T_{*, \ga } (R) \subseteq T$ since $\pi (T_*) \in \Den_*(\bR , 0)$, by Theorem \ref{L11Jan19}.(1). Since 
$$\pi (T) = S_*(\bR ) \in \Den_*(\bR , 0)\;\; {\rm and}\;\; 
 \gb := \ass_R(T) \neq R,$$
 we have that  $T\in \mL_*(R, \gb )$ and $\gb \subseteq \ga$. Since $T_*\subseteq T$,  we have that $\ga \subseteq \gb$, and so $\ga = \gb$. Therefore, $$T\in \mLsRa.$$
  Since $T_*\subseteq T$, we must have $T_*=T$, by the maximality of $T_*$. 

2. Statement 2 follows from statement 1.  

3. The map $\s$ is the composition of the homomorphisms $R\stackrel{\pi}{\ra}\bR\stackrel{\overline{\sigma}}{\ra} \Q_{*, \ga } (R) \simeq Q_*(\bR )$ where $\overline{\sigma}(r+\ga )=\frac{r}{1}$. By \cite[Theorem 2.8.(1)]{larglquot},   $S_*(\bR )=\overline{\sigma}^{-1}(Q_*(\bR )^\times)$. Now, statement 3 follows from statement 1. 
$\Box $\\

{\bf The $\mL_*$-radical $\mL_*\rad (R)$.} 

{\it Definition.} An element $S$ of the set	$\max  \mLsR$ where $*\in \{ l,r,\emptyset\}$ is called the {\em maximal left/right localizable set} and the {\em maximal localizable set}, resp., and the rings $S^{-1}R$, $RS^{-1}$ and $S^{-1}R\simeq RS^{-1}$ are called the {\em maximal left/right localization} and the {\em maximal localization} of $R$, resp..\\

The intersection
\begin{equation}\label{LLradR}
\mL_*\rad (R) = \bigcap_{S\in \max  \mLsR} \ass _R(S)
\end{equation}
is called the $\mL_*$-{\em radical} of $R$. 

For a ring $R$, there is the canonical exact sequence where $ *\in \{ l,r,\emptyset\}$,

\begin{equation}\label{LLradR1}
0\ra \mL_*\rad (R)\ra R\stackrel{\s }{\ra} \prod_{S\in \max \mLsR}S^{-1}R, \;\; \s := \prod_{S\in\max \mLsR}\, \s_S,
\end{equation}
where $\s_S:R\ra S^{-1}R$, $r\mapsto \frac{r}{1}$. A similar sequence exists for $*=$r. \\

{\it Definition.} The sets $\CL \mLsR:=\bigcup_{S\in \mLsR}S$ and $\CN \CL\mLsR := R\backslash \CL \mLsR$ are called the set of $\mL_*$-{\em localizable} and $\mL_*$-{\em non-localizable} elements of $R$, resp., and the intersection
$$ \CC\mLsR =\bigcap_{S\in \max \mLsR}S$$ is called the {\em set of completely * localizable elements} of the ring $R$. \\

By the very definition the sets $\CL \mLsR$, $\CN \CL\mLsR$ and $ \CC\mLsR$ are invariant under the action of the automorphism  group of the ring $R$, i.e. they are {\em characteristic sets}.  \\





{\bf The sets ${}'\mLlR $, $\mLrpR$ and ${}'\mL_{l,r}'(R)$.} For an ideal $\ga $ of a ring $R$, let  
$$ {}'\Den_l(R, \ga) :=\{ S\in \Den_l(R, \ga ) \, | \, S\subseteq {}'\CC_R\}\;\; {\rm and} \;\; \Den_r'(R, \ga) :=\{ S\in \Den_r(R, \ga ) \, | \, S\subseteq \CC_R'\}.$$
 Recall that a  localizable set $S\in \mLsR$ is called {\em perfect} if $\ass_R(S) = \ga (S)$. So, the perfect localizable sets $S$ of $R$ have the `smallest' possible ideal $\ass_R(S)$. 
Proposition \ref{A15Jan19} shows that the sets ${}'\mLlR $, $\mLrpR$ and ${}'\mL_{l,r}'(R)$ are perfect localizable sets  of the ring $R$. See the Introduction for their definitions. 

\begin{proposition}\label{A15Jan19}
We keep the notation as above. 
\begin{enumerate}
\item ${}'\mLlR \subseteq \mL_l^p(R)$ and for each ${}'S\in {}'\mLlR$, $\ass_R(S)= \ga (S)= {}'\pi^{-1}(\ass_l({}'S))$ and  $\RSm \simeq {}'S^{-1}{}'R$ where $\ass_l({}'S)=\{ {}'r\in {}'R\, | \, {}'s{}'r=0$ for some ${}'s\in {}'S\}$. Furthermore, ${}'S\in {}'\Den_l({}'R, \ga / {}'\ga )$ where $\ga = \ass_R(S)$. 
\item $\mLrpR \subseteq \mL_r^p(R)$ and for each $S'\in \mLrpR$, $\ass_R(S)= \ga (S)= \pi'^{-1}(\ass_r(S'))$ and  $\RSm \simeq R'S'^{-1}$ where $\ass_r(S')=\{ r'\in R'\, | \, r's'=0$ for some $s'\in S'\}$. Furthermore, $S'\in \Den_r'(R', \ga / \ga' )$ where $\ga = \ass_R(S)$. 
\item ${}'\mL_{l,r}'(R) \subseteq \mL^p(R)$ and for each $S\in {}'\mL_{l,r}'(R)$, $\ass_R(S)= \ga (S)= {}'\pi^{-1}(\ass_l({}'S))= \pi'^{-1}(\ass_r(S'))$ and  $\RSm \simeq {}'S^{-1}{}'R\simeq R'S'^{-1}$. Furthermore, ${}'S\in {}'\Den_l({}'R, \ga / {}'\ga )$ and $S'\in \Den_r'(R', \ga / \ga' )$  where $\ga = \ass_R(S)$. 
\end{enumerate}
\end{proposition}

{\it Proof}. 1. Recall that  $\ga (S)$ is  the ideal in Proposition \ref{A12Jan19}.(1).

(i) ${}'\ga\subseteq \ga (S)$: The inclusion follows from Proposition \ref{A12Jan19}.(2) (by the minimality of the the ideal ${}'\ga$ since $(S+\ga (S))/\ga (S) \subseteq \CC_{R/ \ga (S)}\subseteq {}'\CC_{R/ \ga (S)}$). 

(ii) $\gb :={}'\pi^{-1}(\ass_l({}'S))\subseteq \ga (S)$: The inclusion $ \gb \subseteq\ga (S)$ follows from the inclusion $ {}'\ga \subseteq \ga (S)$ and the definition of the ideals $\gb$ and $\ga (S)$.

(iii) $\gb = \ass_R(S)$: By Lemma \ref{b15Jan19}, $\ga (S)\subseteq \ass_R(S)$. Now, $ \gb \subseteq \ass_R(S)$, by the statement (ii). Since $\widetilde{S}:=(S+\gb ) / \gb \in \Den_l(R/\gb , 0)$, we must have $\gb \supseteq \ass_R(S)$, by the minimality of the ideal $\ass_R(S)$ (Theorem \ref{L11Jan19}.(3)), and the statement (iii) follows.

 By the statements (ii) and (iii), $\ass_R(S) =\ga (S) = \gb$ (since $\gb\subseteq \ga (S) \subseteq \ass_R(S)=\gb$) and $\RSm\simeq {}'S^{-1}{}'R$. Clearly, ${}'S\in \Den_l({}'R, \ga / {}'\ga )$.

2. Statement 2 can be proven in a similar/dual way as statement 1. 

3. Statement 3 follows from statements  1 and 2.  $\Box $\\

For $L\in \{ {}'\mLlR, \mLrpR, {}'\mL_{l,r}'(R)\}$, let $\ass (L) := \{ \ass_R(S)\, | \, S\in L\}$. Then 

\begin{equation}\label{pLLR}
{}'\mLlR =\coprod_{\ga\in \ass\,  {}'\mLlR }{}'\mLlRa \;\; {\rm where}\;\; 
{}'\mLlRa =\{ S\in {}'\mLlR \, | \, \ass_R(S)=\ga\}, 
\end{equation}

\begin{equation}\label{pLLR1}
 \mLrpR =\coprod_{\ga\in \ass \,  \mLrpR} \mLrpRa  \;\; {\rm where}\;\;\mLrpRa =\{ S\in \mLrpR  \, | \, \ass_R(S)=\ga\},
\end{equation}

\begin{equation}\label{pLLR2}
 {}'\mL_{l,r}'(R) =\coprod_{\ga\in \ass \, {}'\mL_{l,r}'(R) }{}'\mL_{l,r}'(R, \ga )  \;\; {\rm where}\;\; {}'\mL_{l,r}'(R, \ga )=\{ S\in {}'\mL_{l,r}'(R) \, | \, \ass_R(S)=\ga\}.
\end{equation}

{\bf Left/right localizable Ore sets.} Recall that for a ring $R$,  the sets 
$$\LOrelR :=\mLlR \cap \Ore_l(R), \;\;  \LOrerR :=\mLrR \cap \Ore_r(R), \;\;   \LOreR := \mLR \cap \Ore (R)= \Ore (R)$$
are called left, right and Ore localizable, respectively.  Clearly, $\LOreR = \LOrelR \cap \LOrerR$ since 
\begin{eqnarray*}
\LOreR &=&\mLR \cap \Ore (R)=(\mLlR \cap \mLrR ) \cap (\Ore_l (R)\cap \Ore_r(R))\\
& =& (\mLlR \cap \Ore_l (R) ) \cap (\mLrR\cap \Ore_r(R))= \LOrelR \cap \LOrerR .
\end{eqnarray*}
For each elememt $*\in \{ l,r,\emptyset\}$, let $\ass\,  \LOresR :=\{ \ass_R(S)\, | \, S\in \LOresR \}$. Then 

\begin{equation}\label{LOreU}
\LOresR =\coprod_{\ga\in \ass\,  \LOresR }\LOresRa \;\; {\rm where}\;\; 
\LOresRa :=\{ S\in\LOresR  \, | \, \ass_R(S)=\ga\}. 
\end{equation}
Clearly, $\LOresRa =  \mLsRa  \cap \Ore_*(R)$.  Since $\Ore (R)\subseteq \mLR$ see (\ref{LOre=Ore}), we have that $\mL\Ore (R) = \Ore (R)$. So, the letter `$\mL$' is redundant in the definition of the sets `$\mL\Ore (R)$', `$\ass\, \mL\Ore (R)$' and `$\max \mL\Ore (R, \ga )$, and we drop it. So, (\ref{LOreU}) takes the form

\begin{equation}\label{OreU}
\Ore (R) =\coprod_{\ga\in \ass\,  \Ore (R) }\Ore (R,\ga ) \;\; {\rm where}\;\; 
\Ore (R,\ga ):=\{ S\in\Ore (R)  \, | \, \ass_R(S)=\ga\}. 
\end{equation}

\begin{proposition}\label{B15Jan19}

\begin{enumerate}
\item Given $S\in \LOrelR$. Then ${}'\ga \neq R$ (see Proposition \ref{A12Jan19}.(2)) and the ideal ${}'\ga$ is the least ideal $\gb$ of the ring $R$ such that $(S+\gb ) /\gb \in {}'\Den_l(R/ \gb )$. In particular, $(S+ {}'\ga ) / {}'\ga \in {}'\Den_l(R/  {}'\ga )$.
\item  Given $S\in \LOrerR$. Then $\ga' \neq R$ (see Proposition \ref{A12Jan19}.(3)) and the ideal $\ga'$ is the least ideal $\gb$ of the ring $R$ such that $(S+\gb ) /\gb \in {}\Den_r'(R/ \gb )$. In particular, $(S+ \ga' ) / \ga' \in \Den_r'(R/ \ga' )$.
\end{enumerate}
\end{proposition}

{\it Proof}. 1. By Proposition \ref{A12Jan19} and Lemma \ref{b15Jan19}, ${}'\ga \subseteq \ga (S) \subseteq \ass_R(S) \neq R$. Then ${}'S:=(S+ {}'\ga ) / {}'\ga \in {}'\Den_l(R/  {}'\ga )$. Since the ideal ${}'\ga$ is the least ideal of the ring $R$ such that $(S+ {}'\ga ) / {}'\ga \in {}'\CC_{R/{}'\ga}$, statement 1 follows.

2. Statement 2 is proven in a dual way to statement 1. $\Box $\\

{\bf Proof of Proposition \ref{C15Jan19}.} 1. By Proposition \ref{B15Jan19}.(1), $\LOrelR\subseteq {}'\mLlR$. By  Proposition \ref{A15Jan19}.(1), $ {}'\mLlR\subseteq \mL_l^p(R)$.  Now,
\begin{eqnarray*}
\LOrelR  &=& {}'\mLlR\cap \LOrelR = {}'\mLlR\cap \mLR\cap \Ore_l(R) = {}'\mLlR\cap \Ore_l(R),\\
\LOrelR &=&\mL_l^p(R)\cap \LOrelR = \mL_l^p(R)\cap \mLR\cap \Ore_l(R) = \mL_l^p(R)\cap \Ore_l(R).
\end{eqnarray*}

2. Statement 2 can be proven in a similar way to statement 1. 

3. Since $\Ore (R)= \LOrelR\cap \LOrerR$, statement 3 follows from statements 1 and 2 and the fact that $\Ore (R)\subseteq \mLR$, see (\ref{LOre=Ore})). $\Box$\\

Using some of the above results we obtain criteria for a left/right Ore set to be a left/right localizable set, Theorem \ref{B16Jan19} and Theorem \ref{C17Jan19}.\\

{\bf Proof of Theorem \ref{B16Jan19}.} $1\,  (\Rightarrow )$  If $S\in \mLlR$ then $\ga \neq R$. Now, the implication follows from the inclusions ${}'\ga \subseteq \ga (S) \subseteq \ga$ (Lemma \ref{b15Jan19}).

2(a) If ${}'\ga\neq R$ then clearly ${}'S\in {}'\Den ({}'R)$.

$1\,  (\Leftarrow )$ Since ${}'\ga \subseteq \ga (S) \subseteq \ga$ (Lemma \ref{b15Jan19}), the implication follows from the statement 2(a) and Lemma \ref{a15Jan19}.(1).

2(b,c) Since ${}'\ga \subseteq \ga$  (Lemma \ref{b15Jan19}), $\RSm = {}'R\langle {}'S^{-1}\rangle$, by Lemma \ref{a15Jan19}.(2). Now, the statements (b) and (c) follow from the inclusion ${}'S\in {}'\Den_l({}'R)$ (the statement (a)).  $\Box$  

\begin{theorem}\label{C17Jan19}
Let $R$ be a ring, $S\in \Ore_r(R)$, and $\ga = \ass_R(S)$. 
 Then 
\begin{enumerate}
\item  $S\in \mLrR$ iff $\ga'\neq R$ where the ideal $\ga'= \ga'(S)$ of $R$ is as in Proposition \ref{A12Jan19}.(3) and (\ref{ala2}). 
\item Suppose that $\ga'\neq R$. Let $\pi':R\ra R';=R/ \ga'$, $r\mapsto r'=r+\ga'$ and $S'=\pi'(S)$. Then 
\begin{enumerate}
\item $S'\in \Den_r(R')$, 
\item $\ga = (\pi')^{-1}(\ass_r(S'))$,  
\item $RS^{-1}\simeq R'S'^{-1}$, an $R$-isomorphism. 
\end{enumerate}
\end{enumerate}
\end{theorem}

{\it Proof}.  The proof of the theorem is dual to the proof of Theorem \ref{B16Jan19}. $\Box $\\

{\bf Description of the ideal $\ass_R(S)$ for $S\in \mLsR$.} For each localizable set $S\in \mLsRa$, Proposition \ref{B19Jan19} describes the ideal $\ga$. 

\begin{proposition}\label{B19Jan19}
Let $R$ be a ring and $S$ be a multiplicative set in $R$. 
\begin{enumerate}
\item Let  $S\in \mLlRa$. For each pair of elements $s\in S$ and $r\in R$ fix a pair of elements $s_1\in S$ and $r_1\in R$ such that $s_1r-r_1s\in \ga$, and let $\gb_l$ be the ideal of $R$ generated by the elements $s_1r-r_1s$, $\pi_l:R\ra R_l:=R/ \gb_l$, $r\mapsto r+\gb_l$. Then $S_l:=\pi_l(S)\in  \mL\Ore_l(R_l, \ga / \gb_l)$ and $\ga = \pi_l^{-1} (\ga (S_l))$ where the ideal $\ga (S_l)$ of the ring $R_l$ is defined in Proposition \ref{A12Jan19}.(1). 
\item  Let  $S\in \mLrRa$. For each pair of elements $s\in S$ and $r\in R$ fix a pair of elements $s_1\in S$ and $r_1\in R$ such that $rs_1-sr_1\in \ga$, and let $\gb_r$ be the ideal of $R$ generated by the elements $rs_1-sr_1$, $\pi_r:R\ra R_r:=R/ \gb_r$, $r\mapsto r+\gb_r$. Then $S_r:=\pi_r(S)\in  \mL\Ore_r(R_r, \ga / \gb_r )$ and $\ga = \pi_r^{-1} (\ga (S_r))$. 
\item  Let  $S\in \mLRa$, $\gb =\gb_l+\gb_r$ and $\tpi :R\ra \widetilde{R}$, $r\mapsto r+\gb$. Then $\widetilde{S}:=\tpi (S)\in \Ore (\widetilde{R}, \ga / \gb )$ and $\ga = \tpi^{-1} (\ga (\widetilde{S}))$. 
\end{enumerate}
\end{proposition}

{\it Proof}. 1. By the very definition, $\gb_l\subseteq \ga$. By Lemma \ref{a15Jan19}.(2), $S_l\in \mL_l(R_l, \ga / \gb_l)$. By the definition of $\gb_l$, $S_l\in \mL\Ore_l(R_l, \ga / \gb_l)$. By Proposition \ref{C15Jan19}.(1), $\ga / \gb_l= \ga (S_l)$, and so $\ga = \pi_l^{-1} (\ga (S_l))$. 

2 and 3. Statements 2 and 3 can be proven in a similar way. $\Box $\\


Proposition \ref{A19Jan19} gives a sufficient condition for an epimorphic image of a left localizable set to be a left localizable left Ore set. 

\begin{proposition}\label{A19Jan19}
Let $R$ be a ring, $S\in \mL_l(R, \ga )$, $\gb$ be an ideal of $R$ such that $\ga \subseteq \gb$ and $\bR = R/\gb$. Suppose that the left ideal $S^{-1}R\gb$ is an ideal of the ring $S^{-1}R$ such that $S^{-1}R\gb\neq S^{-1}R$. Then $\bS :=(S+\gb ) / \gb\in \mL\Ore_l(\bR ,\overline{\gc }) $ for some ideal $\gc$ of $R$ such that $\gb\subseteq \gc$. 
\end{proposition}

{\it Proof}. Since $S\in \mLlRa$ and $\ga \subseteq \gb$, $\bS \in \Ore_l(\bR)$. Let $\gb_1$ be the kernel of the ring homomorphism $R\ra S^{-1}R / S^{-1}R\gb$, $r\mapsto \frac{r}{1}+S^{-1}R\gb$. Clearly, $\gb \subseteq \gb_1$, $\widetilde{S}=(S+\gb_1)/\gb_1\in \Den_l(R/ \gb_1, 0)$ and $S^{-1}R/S^{-1}R\gb\simeq \widetilde{S}^{-1}(R/\gb_1)$. By Theorem \ref{L11Jan19}.(3), $\bS \in \mL\Ore_l(\bR , \overline{\gc }) $ for some ideal $\gc$ of $R$ such that $\gb\subseteq \gc$. $\Box $


\section{Localization of a ring at an Ore set}\label{LOCORE}

The localization of a ring at an Ore set was introduced in  \cite{larglquot}. The aim of this section is to prove Theorem \ref{10Jan19} which among other things explains what it the localization of a  ring at an Ore set and why it always exists. As a corollary of Theorem \ref{10Jan19} we show that every Ore set is localizable (Theorem \ref{16Jan19}.(1)) and the localization of a ring at an Ore set $S$ is the same as the localization of the ring at the localizable set $S$ (Theorem \ref{16Jan19}.(2)). Proofs of Corollary \ref{c10Jan19} and Proposition \ref{A13Jan19} are given. \\  

{\bf Proof of Theorem \ref{10Jan19}.} 1.(i) $\ga + \ga \subseteq \ga$: Given elements $a,a'\in \ga$. Then $sat=s'a't'=0$ for some elements $s,s',t,t'\in S$. Fix elements $s_1, t_1\in S$ such that $s_1s=\alpha s'$ and $ tt_1=t'\beta$ for some elements $\alpha , \beta \in R$. Then $s_1s, tt_1\in S$ and $$s_1s(a+a')tt_1=s_1(sat)t_1+\alpha (s'a't')\beta=0+0=0,$$ and so $a+a'\in \ga$. 

(ii) $R\ga R\subseteq \ga$: Given elements $r,r'\in R$ and $a\in \ga$. Then $sas'=0$ for some elements $s,s'\in S$. Using the left and right Ore conditions, we have the equalities
$$ s_1r= r_1s\;\; {\rm and}\;\; r_1's_1'= s' r_1'$$ for some elements $s_1, s_1'\in S$ and $r_1, r_1'\in R$. Then 
$$ s_1rar's_1'=r_1sas'r_1'=r_10r_1'=0,$$
and so $rar'\in \ga$, and the statement (ii) follows.

The statements (i) and (ii) imply that the set $\ga$ is an ideal of the ring $R$. 

(iii) $\ga\neq R$: If $\ga = R$ then $1\in \ga$, and so $s\cdot 1 \cdot t=0$ for some elements $s,t\in S$. Then $0=st\in S$, a contradiction. 

2. (i) $S\cap \ga =\emptyset$: If $r\in S\cap \ga$ then $srt=0$ for some elements $s,t\in S$. This is not possible since $0=srt\in S$. 

(ii) $\bS\in \Den (\bR , 0)$: Clearly, $\bS\in \Ore (\bR )$. If $\bs \br =0$ or $\overline{r'}\, \overline{s'}=0$ for some elements $s,s'\in S$ and $r,r'\in R$. Then $sr\in \ga$ and $r's'\in \ga$. Then $s_1srs_1'=0$ or $s_1'r's's_2'=0$ for some elements $s_1, s_2, s_1', s_2'\in S$, and so $r\in \ga$ and $r'\in \ga$, i.e., $\br =0$ and $\overline{r'}=0$, and the statement (ii) follows.  

(iii) $S^{-1}R\simeq \bS^{-1}\bR$ {\em and}  $\ga = \ass (R)$: The statement (iii) follows from Theorem \ref{L11Jan19}.(2,3).

(iv) $\ga = \ga (S) = \ass (R)$: The result follows from the inclusions $\ga \subseteq  \ga (S) \subseteq  \ass (R)$ and 
the equality $\ga = \ass (R)$ (the statement (iii)).

3. (i) $\ga\subseteq \gb$: Given an element $a\in \ga$. Then $sat =0$ for some elements $s,t\in S$, and so $s^\dag a^\dag t^\dag =0$. Hence $a^\dag =0$ in $R^\dag $ since $s^\dag ,t^\dag \in \CC_{R^\dag }$, and so $a\in \ga$. 

(ii) {\em The map $\bS^{-1} \bR \ra {S^\dag}^{-1} R^\dag $, $\bs^{-1}\br\mapsto {s^\dag}^{-1}r^\dag $ is a ring epimorphism:} The statement (ii) follows from the universal property of localization.

4(a) (i)  $Q= \{ f(s)^{-1} f(r) \, | \, s\in S, r\in R\} =  \{  f(r)f(s)^{-1} \, | \, s\in S, r\in R\}$: The statement (i) follows from the fact that $S\in \Ore (R)$ and that the ring $Q$ is generated by the sets $f(R)$ and $\{ f(s)^{-1} \, | \, s\in S\}$. 

(ii) $\ga \subseteq \ker (f)$: Given an element $a\in \ga$. Then $sat=0$ for some elements $s,t\in S$, and so 
$$ 0= f(sat) = f(a) f(s) f(t).$$
Hence, $f(a)=0$ since $f(s), f(t) \in Q^\times$, and so $a\in \ker (f)$, and the statement (ii) follows. 

(iii) {\em The map $\bS^{-1}\bR \ra Q$, $ \bs^{-1} \br \mapsto f(s)^{-1} f(r)$ is a ring epimorphism:} Since the ring $Q$ is generated by the set $f(R)$ and $\{ f(s)^{-1} \, | \, s\in S\}$, and $\ga \subseteq \ker (f)$, the statement (iii) follows from the universal property of localization.

(iv) {\em The kernel of the map in the statement (iii) is} $\bS^{-1} ( \ker (f) / \ga )$: The statement is obvious.

4(b) The statement (b) follows from the statement (a). $\Box$\\

{\em Definition, \cite{larglquot}.} Let $S\in \Ore (R)$. The ring $\bS^{-1}\bR \simeq \bR\bS^{-1}$ in Theorem \ref{10Jan19}.(2) is called the {\em localization of the ring $R$ at the Ore set} $S$ and is denoted by $S^{-1}R=RS^{-1}$. The ideal $\ga$ in Theorem \ref{10Jan19}.(1) is denoted by $\ass (S)$. \\

Statement 4 (and statement 3) of Theorem \ref{10Jan19} is the universal property of localization of a ring at Ore set. If the Ore set is a denominator set this is precisely the universal property of localization at a denominator set. 

As a corollary of Theorem \ref{10Jan19} we have that every Ore set is a localizable set (Theorem \ref{16Jan19}.(1)). This result was proven in \cite{larglquot} by a different method. 

\begin{theorem}\label{16Jan19}
Let $R$ be a ring. 
\begin{enumerate}
\item $\Ore (R)\subseteq \mLR$, i.e., every Ore set of $R$ is a localizable set of $R$. 
\item If $S\in \Ore (R)$ then $S\in \mLRa$ where the ideal $\ga$ is the ideal in Theorem \ref{10Jan19}.(1) and the localization of the ring $R$ at the Ore set $S$ coincides with the localization of the ring $R$ at the localizable set. 
\end{enumerate}
\end{theorem}

{\it Proof}. 1. Statement 1 follows from statement 2. 

2. We keep the notation of Theorem \ref{10Jan19}. Let $\ga$ be the ideal in Theorem \ref{10Jan19}.(1).  By Theorem \ref{10Jan19}.(2), $\bS\in \Den (\bR , 0)$ where $\bS = (S+\ga ) / \ga $ and $\bR = R/\ga$. Therefore, $\ga (S)=\ga$ where the ideal $\ga (A)$ is as in Proposition \ref{A12Jan19}.(1) (the inclusion $\ga \subseteq \ga (S)$ follows from Theorem \ref{10Jan19}.(1), and the inclusion $\ga (S) \subseteq \ga$ follows from the fact that $\bS\subseteq \CC_{\bR }$). Since $ \ga = \ga (S) \subseteq \ass_R(S)$ (Lemma \ref{b15Jan19}) and $\bS \in \Den (\bR , 0)$, we have that $\ga = \ass_R(S)$, $S\in \mLRa$ and the localization of $R$ at the Ore set $S$ coincides with the localization of $R$ at the localizable set $S$.  $\Box $\\

{\bf Proof of Corollary \ref{c10Jan19}.} 1. Statement 1 follows from Theorem \ref{10Jan19}.(1).

2. (i) $S_l\in \Ore_l(R_l)$ {\em with} $ \ass_l(S_l)=0$: Since $S\in \Ore (R)$ and $R_l=R/\ass_l(S)$, the statement (i) is obvious.

(ii) $S_l\in \Den_r'(R_l, \ga / \ga_l)$: Since $S\in \Ore (R)$, we have that $S_l\in \Den_r(R_l)$ since $\ass_l(S_l)=0\subseteq \ass_r(S_l)$, by the statement (i). By Theorem \ref{10Jan19}.(1), $\ass_r(S_l) = \ga / \ga_l$. In more detail, $\pi_l(r)\in \ass_r(S_l)$ iff $0=\pi_l(r) \pi_l(t) = \pi_l(rt)$ for some element $t\in S$ iff $rt\in \ga_l$ iff $srt=0$ for some $s\in S$ iff $r\in \ga$, by Theorem \ref{10Jan19}.(1). 

(iii) $R_lS_l^{-1}\simeq S^{-1}R\simeq RS^{-1}$: The statement (iii) follows from Theorem \ref{10Jan19}.(2). 

3. Statement 3 is proven in a similar fashion as statement 2. 

4. Let ${}'\ga = {}'\ga (S)$ and $\ga'=\ga'(S)$. By the very definition of the ideals $\ga_l$ and $\ga_r$, ${}'\ga \supseteq \ga_r$ and $\ga'\supseteq \ga_l$. The reverse inclusions follow from statements 2 and 3 and the `minimality of the ideals ${}'\ga $ and $\ga'$ in the sense of Lemma \ref{A12Jan19}.(2,3), resp.. $\Box$\\ 

{\bf Proof of Proposition \ref{A13Jan19}.} (i) $\widetilde{S}$ {\em is a multiplicative set of $R$ such that} $S\subseteq \widetilde{S}$: The set $\bS = \pi (S)$ is a multiplication set of the ring $\bR$, hence the set $\widetilde{S}=\pi^{-1}(\bS )$ is a multiplicative set of $R$. 

(ii) $\widetilde{S}\in \mLR$: It suffices to show that $\widetilde{S}\in \mLlR$ since then by symmetry we will get $\widetilde{S}\in \mLrR$. Given elements $s\in \widetilde{S}$ and $r\in R$. Since $\pi (\tilde{s})= \bS \in \Den_l(\bR , 0)$, $ \bs_1\br = \br_1\bs$ for some elements $\bs_1\in \bS$ and $ \br_1\in \bR$. Then $s_1r-r_1s\in \ga$. 
 By Theorem \ref{10Jan19}.(1),  $\ga \subseteq \ass_R(\widetilde{S})$,  $\widetilde{S}\in \mLlR$. 
 

(iii) $\widetilde{S}\in \mLRa$: Let $\tilde{\ga}=\ass_R(\widetilde{S})$. We have to show that $\tilde{\ga}=\ga$. Since $S\subseteq \widetilde{S}$, we have the inclusion $\ga \subseteq \tilde{\ga }$ (Lemma \ref{a12Jan19}). Since $\pi (\widetilde{S})=\bS \in \Den_l(R, 0)$, we have the reverse inclusion $\ga \supseteq \tilde{\ga }$ (by Theorem \ref{10Jan19}.(4) and Theorem \ref{16Jan19}), and so $\tilde{\ga}=\ga$.

(iv) $\widetilde{S}^{-1}R\simeq S^{-1}R$ {\em is an $R$-isomorphism}: By the statement (iii), $\widetilde{S}^{-1}R\simeq \pi(S)^{-1}\bR = \bS^{-1}\bR = S^{-1}R$. $\Box$ \\

The following obvious lemma is a useful criterion when an epimorphism image of Ore set is an Ore set. 

\begin{lemma}\label{a19Jan19}
Let $R$ be a ring, $S\in \Ore (R)$ and $\gb$ be an ideal of $R$. Then $\bS :=(S+\gb ) / \gb \in \Ore (R/ \gb )$ iff $S\cap \gb = \emptyset$. 
\end{lemma}

{\it Proof}. Straightforward. $\Box $

An element $a$ of ring $R$ is called a {\em normal element} if $aR=Ra$. We keep the notation of Corollary \ref{c10Jan19}. 
 The next example shows that for an Ore set $S$ of $R$, the ideals $\ga_l$, $\ga_r$ and $\ga$ are distinct and $\ga = \ga_l+\ga_r$, but, in general,  $\ga \varsupsetneqq  \ga_l+\ga_r$ (see the second example).

{\em Example}. Let $P= K[x_1, x_2, \ldots ,  y_1, y_2, \ldots ]$ be a polynomial algebra over a field $K$ in countably many variables $x_1, x_2, \ldots ,  y_1, y_2, \ldots $. Let $R$ be a $K$-algebra generated by $P$ and an element $a$ subject to  the defining  relations:
$$ ax_1=0, \;\; ax_i=x_{i-1}a\;\; (i\geq 2), \;\; y_1a=0, \;\; y_ia= ay_{i-1}\;\; (i\geq 2).$$
Then $a$ is a normal element of $R$ (since $aR=\sum_{i\geq 1}P'a^i=Ra$ where $P'=K[x_1, \ldots , y_2, \ldots ]$), hence $S_a:=\{ a^i\, | \, i\in \N\}\in \Ore (R)$, $\ass_l(S_a)= (x_1, x_2, \ldots )$, $\ass_r(R)= (y_1, y_2, \ldots )$, $\ass(S_a)=(x_1,  \ldots , y_1,  \ldots )= \ass_l(S_a)+\ass_r(S_a)$ and $S_a^{-1}R\simeq RS_a^{-1}\simeq K[a,a^{-1}]$ (Corollary \ref{c10Jan19}.(2,3)). Clearly, the ideals $\ass_l(S_a)$, $\ass_r(S_a)$ and $\ass(S_a)$ are distinct and  $\ass(S_a)= \ass_l(S_a)+\ass_r(S_a)$. The next example shows that, in general, $\ass (S_a)$ properly holds the sum  $\ass_l(S_a)+\ass_r(S_a)$.

{\em Example}. Let $P= K[x_1, x_2, \ldots , \ldots , y_1, y_2, \ldots ]$ be a polynomial algebra over a field $K$ in variables $x_1, x_2, \ldots ,  y_1, y_2, \ldots $. Let $R$ be a $K$-algebra generated by $P$ and an element $a$ subject to  the defining relations:
$$ ax_1=0, \;\; ax_i=x_{i-1}a\;\; (i\geq 2), \;\; y_1a=0, \;\; y_2a=a(y_1+x_2), \;\; y_ia= ay_{i-1}\;\; (i\geq 3).$$
Then $a$ is a normal element of $R$ (since $aR\subseteq Ra$ (as $ay_1=(y_2-x_1)a$) and $Ra\subseteq aR$), hence $S_a:=\{ a^i\, | \, i\in \N\}\in \Ore (R)$, $\ass_l(S_a)= (x_1, x_2, \ldots )$, $\ass_r(R)= (y_0)$, $\ass(S_a)=(x_1, x_2, \ldots ,  y_1, y_2, \ldots )\varsupsetneqq \ass_l(S_a)+\ass_r(S_a)$ and $S_a^{-1}R\simeq RS_a^{-1}\simeq K[a,a^{-1}]$ (Corollary \ref{c10Jan19}.(2,3)).


\section{Localization of a ring at an almost Ore set}\label{LAORE}

The aim of this section is to introduce almost Ore sets, to give a criterion for an almost Ore set to be a localizable set (Theorem \ref{A16Jan19}.(1), 
 Theorem \ref{D16Jan19}.(2) and Theorem \ref{E16Jan19}.(2)) and for each localizable almost Ore set $S$ to give an explicit description of the ideal $\ass_R(S)$ and of the ring $S^{-1}R$ (Theorem \ref{A16Jan19}.(2),  Theorem \ref{D16Jan19}.(3) and Theorem \ref{E16Jan19}.(3)). \\

{\bf Almost Ore sets.} Let $R$ be a ring. Recall that a multiplicative set $S$ of $R$ is a {\em left} (resp., {\em right}) {\em Ore set} if the left (resp., right) Ore condition holds: For any elements $s\in S$ and $r\in R$, there are elements $s_1\in S$ and $r_1\in R$ such that $s_1r=r_1s$ (resp., $rs_1=sr_1$). A left and right Ore set of $R$ is called an {\em Ore set} of $R$. \\

{\it Definition.} Let $R$ be a ring. A multiplication set $S$ of $R$ is called an {\em almost left } (resp., {\em right}) {\em Ore set} of $R$ if the {\em almost left} (resp., {\em right}) {\em Ore condition} holds: \\

(ALO) For any elements $s\in S$ and $r\in R$ there are elements $s_1, s_2\in S$ and $r_1\in R$ such that $(s_1r-r_1s) s_2=0$. \\

(ARO)  For any elements $s\in S$ and $r\in R$ there are elements $s_1, s_2\in S$ and $r_1\in R$ such that $s_2(rs_1-sr_1) =0$. \\

Let $\AOrelR$ and $\AOrerR$ be the sets of almost left and right Ore sets, respectively. Elements of the set ${\rm AOre}(R) := \AOrelR\cap \AOrerR $ are called {\em almost Ore sets}.

Clearly, $\Ore_*(R)\subseteq \AOresR$ for $*\in \{ l,r,\emptyset\}$. Let 
\begin{eqnarray*}
\LAOlR &:=& \mLlR\cap \AOrelR, \\
\LAOrR &:=& \mLrR\cap \AOrerR,\\
\LAOR &:=& \LAOlR\cap  \LAOrR= \mLR\cap \AOreR .
\end{eqnarray*}


Elements of the sets $\LAOlR$ and $\LAOrR$ are called {\em left localizable almost left Ore sets} and {\em right  localizable almost  right Ore sets}, respectively. Elements of the set $\LAOR$  are called {\em localizable almost  Ore sets}.

Let $\ass \, \LAOsR :=\{ \ass_R(S) \, | \, S\in \LAOsR\}$. Then 
\begin{equation}\label{LAOaU}
\LAOsR =\coprod_{\ga \in \ass \, \LAOsR} \mL{\rm AO}_*(R, \ga )\;\; {\rm where}\;\; \mL{\rm AO}_*(R, \ga ):=\{ S\in \LAOsR\, | \, \ass_R(S) = \ga \}. 
\end{equation}
Clearly, $ \mL{\rm AO}_*(R, \ga )= \LAOsR\cap \mLsRa$. \\

Theorem \ref{A16Jan19}.(1) is a criterion for an almost Ore set to be a localizable set. It also gives an explicit description of the ideal $\ass_R(\ga )$ for each $S\in \mL{\rm AO}(R)$ (Theorem \ref{A16Jan19}.(2b)). 

\begin{theorem}\label{A16Jan19}
Let $R$ be a ring, $S\in \AOreR$, $\ga_l$ and $\ga_r$ be the ideals of $R$ generated by the sets $\ass_l(S)=\{ r\in R\, | \, sr=0$ for some $s\in S\}$ and $\ass_r(S)=\{ r\in R\, | \, rs=0$ for some $s\in S\}$, respectively. 
\begin{enumerate}
\item The following statements are equivalent.
\begin{enumerate}
\item $S\in \mLR$. 
\item $\ga_l+\ga_r\neq R$. 
\end{enumerate}
\item Suppose that $\tilde{\ga} = \ga_l+\ga_r\neq R$. Let $\widetilde{R}=R/\tilde{\ga}$ and $ \widetilde{S} = (S+\tilde{\ga})/\tilde{\ga}$. Then 
\begin{enumerate}
\item $\widetilde{S}\in \Ore (\widetilde{R})$. 
\item $\ass_R(S)=\tilde{\pi}^{-1} (\ga^\circ)$ where $\tilde{\pi}:R\ra \widetilde{R}$, $r\mapsto \tilde{r}=r+\tilde{\ga}$ and $\ga^\circ$ is the ideal in Theorem \ref{10Jan19}.(1) for the Ore set $\widetilde{S}\in \Ore_l(\widetilde{R})$,  that is $\ga^\circ =\{ \tilde{r}\in \widetilde{R}\, | \, \tilde{s}\tilde{r} \tilde{t}=0$ for some elements $\tilde{s}, \tilde{t}\in \widetilde{S}\}$.
\item Let $\ga = \ass_R(S)$ and $\pi :R\ra \bR := R/\ga$, $r\mapsto \br = r+\ga$. Then $\bS :=\pi (S)\in \Den (R, 0)$.
\item $S^{-1}R\simeq \bS^{-1}\bR$, an $R$-isomorphism. 
\end{enumerate}
\end{enumerate}
\end{theorem}

{\it Proof}. Recall that $\ga = \ass_R(S)$ and $\ga (S)$ be the least ideal of the ring $R$ such that $(S+\ga (S))/ \ga (S) \subseteq \CC_{R/\ga (S)}$ (Proposition \ref{A12Jan19}.(1)). 

(i) $\tilde{\ga} \subseteq \ga (S) \subseteq \ga$: The first inclusion follows from (\ref{ala}) and the second one does from Lemma \ref{b15Jan19}.

$1\, (a\Rightarrow b)$ If $S\in \mLR$ then $\ga \neq R$, and so $\tilde{\ga } \neq R$, by the statement (i). 

$1\, (b\Rightarrow a)$ 

(ii) {\em If $\tilde{\ga } \neq R$ then} $\widetilde{S}\in \Ore (\widetilde{R})$: The statement (ii) follows at once from the conditions (ALO) and (ARO).

Since each Ore set is localizable, the implication $ (b\Rightarrow a)$ follows from the statement (ii).

2(a) The statement (a) is the same as the statement (ii).

2(b,c,d) The statements (b) , (c) and (d) follow from the statement (a) and Theorem \ref{10Jan19}.(1,2). $\Box $\\

{\bf Criterion for an almost left/right Ore set to be a left localizable set.} Theorem \ref{D16Jan19} is such a criterion.

\begin{theorem}\label{D16Jan19}
Let $R$ be a ring, $S\in \AOrelR$, $\ga = \ass_R(S)$, $\ga_r$ be the ideal of $R$ generated by the set $\ass_r(S)= \{ r\in R\, | \, rs=0$ for some $s\in S\}$. Then 
\begin{enumerate}
\item $S_r:= (S+\ga_r)/\ga_r\in \Ore_l(R/ \ga_r)$. 
\item $S\in \mLlR$ iff $S_r\in \mL_l(R/ \ga_r)$ iff ${}'\ga (S_r)\neq R/\ga_r$ iff ${}'\ga (S)\neq R$ (see Proposition \ref{A12Jan19}.(2) and (\ref{ala1})).
\item Suppose that ${}'\ga := {}'\ga (S)\neq R$. Let ${}'\pi: R\ra {}'R=R/{}'\ga$, $r\mapsto {}'r= r+{}'\ga$ and ${}'S={}'\pi (S)$. Then 
\begin{enumerate}
\item ${}'S\in {}'\Den_l({}'R)$. 
\item $\ga = {}'\pi^{-1}(\ass_l({}'S))$.
\item $S^{-1}R\simeq {}'S^{-1}{}'R$, an $R$-isomorphism. 
\end{enumerate}
\end{enumerate}
\end{theorem}

{\it Proof}. 1. Statement 1 follows at once from the condition (ALO) and the definition of the ideal $\ga_r$. 

2. Since $\ga_r \subseteq \ga (S) \subseteq \ass_R(S)$ (Lemma \ref{b15Jan19}), $S\in \mLlR$ iff $S_r\in \mL_l(R/\ga_r)$, by Lemma \ref{a15Jan19}.(1). By Theorem \ref{B16Jan19}.(1), $S_r\in \mL_l(R/ \ga_r)$ iff ${}'\ga (S_r) \neq R/\ga_r$. Since $\ga_r\subseteq {}'\ga (S)$, ${}'\ga (S)/ \ga_r= {}'\ga (S_r)$. Therefore, ${}'\ga (S_r)\neq R/ \ga_r$ iff ${}'\ga (S) \neq R$.

3. Since $\ga_r \subseteq {}'\ga (S)\subseteq \ga (S))\subseteq \ass_R(S)$, $\RSm \simeq (R/ \ga_r)\langle S_r^{-1}\rangle$ (Lemma \ref{a15Jan19}.(1)) and ${}'\ga (S)/\ga_r={}'\ga (S_r)$, statement 3 follows from Theorem \ref{B16Jan19}.(2). 
 $\Box $

\begin{theorem}\label{E16Jan19}
Let $R$ be a ring, $S\in \AOrerR$, $\ga = \ass_R(S)$, $\ga_l$ be the ideal of $R$ generated by the set $\ass_l(S)= \{ r\in R\, | \, sr=0$ for some $s\in S\}$. Then 
\begin{enumerate}
\item $S_l:= (S+\ga_l)/\ga_l\in \Ore_r(R/ \ga_l)$. 
\item $S\in \mLrR$ iff $S_l\in \mL_r(R/ \ga_l)$ iff $\ga' (S_l)\neq R/\ga_l$ iff $\ga' (S)\neq R$ (see Proposition \ref{A12Jan19}.(3) and (\ref{ala2})).
\item Suppose that $\ga' := \ga' (S)\neq R$. Let $\pi': R\ra R'=R/\ga'$, $r\mapsto r'= r+\ga'$ and $S'=\pi' (S)$. Then 
\begin{enumerate}
\item $S'\in \Den_r'(R')$. 
\item $\ga = \pi'^{-1}(\ass_r(S'))$.
\item $RS^{-1}\simeq R'S'^{-1}$, an $R$-isomorphism. 
\end{enumerate}
\end{enumerate}
\end{theorem}

{\it Proof}. The proof of the theorem is `dual' to the proof of Theorem \ref{D16Jan19}. $\Box$ \\

{\bf  Almost Ore sets are perfect localizable sets.} 
Corollary below shows that almost Ore sets are perfect localizable sets.

\begin{proposition}\label{X17Jan19}

\begin{enumerate}
\item $\LAOlR \subseteq \pmLlR \subseteq \mL^p_l(R)$ and  $\LAOlR = \pmLlR \cap \AOrelR = \mL^p_l(R)\cap \AOrelR$.
\item $\LAOrR \subseteq \mLrpR \subseteq \mL^p_r(R)$ and  $\LAOrR = \mLrpR \cap \AOrerR= \mL^p_r(R)\cap \AOrerR$.
\item  $\LAOR \subseteq {}'\mL_{l,r}'(R) \subseteq \mL^p (R)$.
\end{enumerate}
\end{proposition}

{\it Proof}. 1. Statement 1 follows from Theorem \ref{D16Jan19}.(3).

2. Statement 2 follows from Theorem \ref{E16Jan19}.(3).

3. Statement 3 follows from statements 1 and 2. $\Box$


\section{ Classification of maximal localizable sets and maximal Ore sets in semiprime Goldie ring}\label{MAXORE}

The aim of this section is to classify the maximal Ore sets in a semiprime Goldie ring (Theorem \ref{17Jan19}.(1)). One of the key results that is used in the proof of Theorem \ref{17Jan19} is Theorem \ref{A17Jan19}. The concept the maximal left denominator set of a ring was introduced and studied in \cite{larglquot}.

\begin{lemma}\label{a17Jan19}
Let $R$ be  a ring and $S\in \Ore (R)$. For any elements $s,t\in S$, there is an element $\nu\in S$ such that $\nu=x_1s=x_2t=sy_1=ty_2$ for some elements $x_i, y_i\in R$ for $i=1,2$. 
\end{lemma}

{\it Proof}. Since $S\in \Ore_l (R)$, $s_1s=r_1t$ for some elements $s_1\in S$ and $r_1\in R$. Since $S\in \Ore_r (R)$, $ss_2=tr_2$ for some elements $s_2\in S$ and $r_2\in R$.  Then $\nu = ss_2s_1s\in S$ and $$ \nu = s\cdot (s_2s_1s)=ss_2s_1\cdot s= t\cdot r_2s_1s= ss_2r_1\cdot t, $$
as required.  $\Box $

\begin{theorem}\label{A17Jan19}
Let $R$ be a semiprime left Goldie ring. Then $\CC_R\cap \ass_R(S)=\emptyset$ for all $S\in \Ore (R)$. 
\end{theorem}

{\it Proof}. The ring $R$ is a semiprime left Goldie ring. By Goldie's Theorem,  its classical left quotient ring $Q= Q_{l, cl} (R):=\CC_R^{-1}R$ is a semisimple Artinian ring, $Q=\prod_{i=1}^n Q_i$ where $Q_i$ are simple Artinian rings.  The map 
$$\s : R\ra Q=\prod_{i=1}^n Q_i, \;\; r\mapsto \frac{r}{1}=(r_1, \ldots , r_n)$$ is a ring monomorphism. We identify the ring $R$  via $\s$ with its image in $Q$. So, $r=\frac{r}{1}=(r_1, \ldots , r_n)$ where $r_i=\frac{r}{1}\in Q_i$. 

Suppose that $\CC_R\cap \ass_R(S)\neq \emptyset$ for some $S\in \Ore (R)$, we seek a contradiction. Fix an element $c$ in the intersection $\CC_R\cap \ass_R(S)$. By Theorem \ref{10Jan19}.(1), 
$$sct=0$$ for some elements $s,t\in S$. By Lemma \ref{a17Jan19}, we can assume that $s=t$, i.e., $scs=0$, i.e., $Q_i\ni s_ic_is_i=0$ for all $i=1, \ldots , n$ where $s= (s_1, \ldots , s_n) $ and $c=(c_1, \ldots , c_n)$. Clearly, $c_i\in Q_i^\times$ for all $i=1, \ldots , n$. Hence, $$s_i\not\in Q_i^\times\;\; {\rm  for\; all}\;\; i=1, \ldots , n$$ 
(since $s_ic_is_i=0$).  The ideal of $R$,  $\ga_l=\ass_l(S)$, is a nonzero ideal (since $s\cdot cs=0$ and $0\neq cs \in \ga_l$).  The ring $Q$ is a left Noetherian ring.  Hence $\CC_R^{-1}\ga_l$ is an ideal of $Q$. If $s_i=0$ then 
$\CC_R^{-1}\ga_l\cap Q_i=Q_i$. Suppose that $s_i\neq 0$. Then  $0\neq c_is_i\in  \CC_R^{-1}\ga_l\cap Q_i$, and so 
 $$0\neq \CC_R^{-1}\ga_l\cap Q_i=Q_i$$
since  $Q_i$ is a simple ring. Therefore, $\CC_R^{-1}\ga_l\cap Q_i=Q_i$ for all $i=1, \ldots , n$, and so 
$$\CC_R^{-1}\ga_l=\prod_{i=1}^n Q_i= Q.$$
 Hence, $\CC_R\cap \ga_l\neq \emptyset$. Fix an element $a\in \CC_R\cap \ga_l$. Then $s'a=0$ for some $s'\in S$ (since $a\in \ga_l= \ass_l(S))$  but $a\in \CC_R$, hence $0=s'\in S$, a contradiction.  $\Box $\\

{\bf Classification of maximal Ore sets of a semiprime Goldie ring.} \\


{\bf Proof of Theorem \ref{17Jan19}}. 1. By \cite[Theorem 4.1]{Bav-Crit-S-Simp-lQuot}, 
$$\max \Den (R)=\{ \CC (\gp ) \, | \, \gp \in \min (R)\}.$$ By Theorem \ref{18Jan19}, $$\{ \CC (\gp ) \, | \, \gp \in \min (R)\}=\CN_*\;\; {\rm for}\;\; *\in \{ l,r,\emptyset\}.$$
To finish the proof of statement 1 it suffices to show that every Ore set $S$ of the ring $R$ is contained in a maximal denominator set. 

The ring $R$ is a semiprime Goldie ring. The left and right quotient ring of $R$, $$Q_{cl}(R)= \CC_R^{-1}R\simeq R\CC_R^{-1}\simeq \prod_{i=1}^n Q_i$$ is a product of simple Artinian ring $Q_i$. Let $\ga = \ass_R(S)$.  By Theorem \ref{A17Jan19}, $$\CC_R^{-1}\ga \neq Q.$$ Hence, up to order, $\CC_R^{-1}\ga = \prod_{i=m+1}^n Q_i$  for some $m$ such that $1\leq m <n$. We have ring homomorphisms
$$ \s : R\stackrel{\pi}{\ra}R/ \ga \stackrel{\tau}{\ra}\CC_R^{-1}(R/ \ga )\simeq \overline{Q}:=\prod_{i=1}^m Q_i $$
where $\s = \tau \pi$, $\pi (r) =\br := r+\ga$ and $\tau (\br ) = \frac{\br}{1}$. The homomorphism $\pi$ is an epimorphism and the homomorphism $\tau$ is a monomorphism. Let $\bS = \pi (S)$. By Theorem \ref{10Jan19}.(2), $\bS\in \Den (R/ \ga , 0)$. In particular, $\bS\subseteq  \overline{Q}^\times$. Therefore, $\s_1(S) \subseteq Q_1^\times$ where $\s_1:R\ra Q\ra Q_1$ where the first map 
 is $r\mapsto \frac{r}{1}$ and the second is the projection onto $Q_1$. Therefore, $$S\subseteq  \s_1^{-1} (Q_1^\times ).$$
  By the explicit description of the set $\max \Den (R)$ at the beginning of the proof, $\s_1^{-1} (Q_1^\times )\in \max \Den (R)$, as required.

2. In view of the first equality in statement 1, statement 2 is  \cite[Theorem 4.1.(2d)]{Bav-Crit-S-Simp-lQuot}.

3. In view of the first equality in statement 1, statement 3 is  \cite[Theorem 4.1]{Bav-Crit-S-Simp-lQuot}.

4. In view of the first equality in statement 1, statement 4 is \cite[Theorem 4.1.(2c)]{Bav-Crit-S-Simp-lQuot}. $\Box $


\section{Localization of a module at a localizable set}\label{LOCMOD}

The aim of the section is to introduce the concept of localization of a module at a localizable set and to consider its basic properties.\\

{\it Definition.} Let $R$ be a ring, $S\in \mLsRa$ where $*\in \{ l , \emptyset \}$ and $M$ be a left  $R$-module. Then $S^{-1}M:= S^{-1}R\t_R M$ is called the {\em localization} of $M$ at $S$. If $S\in \mLsRa$ where $*\in \{ r , \emptyset \}$ and $M$ be a right  $R$-module. Then $MS^{-1}:= M\t_R RS^{-1}$ is called the {\em localization} of $M$ at $S$.  \\

We consider the case when $*\in \{ l , \emptyset \}$ and $M$ is a left  $R$-module. By the very definition,  $S^{-1}M$ is a left  $S^{-1}R$-module. The elements of the $S^{-1}R$-module $S^{-1}M$ are denoted by $s^{-1}m$. In particular, $s^{-1}r\t m=s^{-1}rm$ for $s\in S$, $r\in R$ and $m\in M$, and $\frac{m}{1}:= 1\t m$. The map $$i_M: M\ra S^{-1}M, \;\; m\mapsto 1\t m$$ is an $R$-homomorphism.

\begin{proposition}\label{A20Jan19}
Let $R$ be a ring, $S\in \mLsRa$ where $*\in \{ l , \emptyset \}$, $M$ be an $R$-module, and $i_M : M\ra S^{-1}M$, $m\mapsto 1\t m$. Then 
\begin{enumerate}
\item $S^{-1}M\simeq \bS^{-1}(M/ \ga M)$ where $\bS := (S+\ga ) / \ga \in \Den_l(\bR , 0)$ and $\bR = R/ \ga$ (Theorem \ref{L11Jan19}.(1)). 
\item Let $\CM$ be an $S^{-1}R$-module and $f: M\ra \CM$ be an $R$-homomorphism. Then there is a unique $S^{-1}R$-homomorphism $S^{-1}f: S^{-1}M\ra \CM$ such that $f=S^{-1}f\circ i_M$. 
\item $\gt_S(M):= \ker (i_M)=\{ m\in M \, | \, sm\in \ga M$ for some $s\in S\}$. 
\end{enumerate}
\end{proposition}

{\it Proof}. 1. Let $\bM=M/\ga M$. Then $S^{-1}(\ga M) = S^{-1}R\t_R \ga M = S^{-1}R\ga \t_RM=0$, and so $$S^{-1}M= S^{-1}\bM = S^{-1} R\t_R\bM\simeq \bS^{-1}\bR \t_{\bR } \bM = \bS^{-1} \bM . $$

2. The $R$-homomorphism $f: M\ra \CM$  determines the ring homomorphism $R\ra \End_\Z (\CM )$, $ r\mapsto (m\mapsto rm)$. The images of the elements of the set $S$ in $\End_\Z (\CM )$ are units. Now, statement 2 follows from Theorem \ref{L11Jan19}.(4). 

3. Statement 3 follows from statement 1.  $\Box$\\

For a ring $R$, let $R-{\rm mod}$ be the category of left $R$-modules.  By Proposition \ref{A20Jan19}.(1), the localization at $S$ functor, $$S^{-1} : R-{\rm mod}\ra S^{-1}R-{\rm mod}, \;\; M\mapsto S^{-1}M,$$ is the composition of two right exact functors 
\begin{equation}\label{Slcm}
S^{-1} = \bS^{-1} \circ (R/\ga \t_R-).
\end{equation}
Therefore, the functor $S^{-1}$ is also a {\em right exact functor} for all $S\in \mLsR$: Given a short exact sequence of $R$-modules $0\ra M_1\ra M_2\ra M_3\ra 0$, then the sequence of $S^{-1}R$-modules 

\begin{equation}\label{Slcm1}
0\ra \bS^{-1} (M_1\cap \ga M_2/\ga M_1)\ra S^{-1}M_1\ra S^{-1}M_2\ra  S^{-1}M_3\ra 0
\end{equation}
is exact. Notice that $$M_1\cap \ga M_2/\ga M_1\simeq \bM_1\cap (\ga M_2/\ga M_1).$$ 

An $R$-module $M$ is called $S$-{\em torsion} (resp., $S$-{\em torsionfree}) if $S^{-1}M=0$ (resp., $\gt_S(M)=0$, i.e., $i_M: M\ra S^{-1}M$, $m\mapsto 1\t m$ is an $R$-module monomorphism). Let $\gf_S (M)=\im (i_M)$, the image of the map $i_M$, and we have a short exact sequence of $R$-modules 
\begin{equation}\label{Slcm2}
0\ra \gt_S(M)\ra M\ra \gf_S(M)\ra 0. 
\end{equation}
Clearly, $\ga M \subseteq \gt_S(M)$, 
$$\gt_S(M)/\ga M=\tor_{\bS } (\bM ):=\{ \overline{m}\in \bM \, | \, \bs \overline{m} =0\;\; {\rm  for\;  some\; element} \;\; \bs \in \bS\}$$  where $\bM = M/ \ga M$ and 
$$\gf_S(M) = M/\gt_S(M)\simeq (M/\ga M) / (\gt_S(M)/ \ga M)\simeq \bM / \tor_{\bS } (\bM ).$$
By taking the short exact sequence (\ref{Slcm2}) modulo $\ga M$, we obtain a short exact sequence of $\bR$-modules 
\begin{equation}\label{Slcm3}
0\ra \tor_{\bS}(\bM )\ra \bM\ra \bM/\tor_{\bS}(\bM )\ra 0. 
\end{equation}

\begin{lemma}\label{a20Jan19}
Let $R$ be a ring, $S\in \mLsRa$ where $*\in \{ l, \emptyset\}$ and $M$ be an $R$-module. Then 
\begin{enumerate}
\item $\gt_S\gf_S(M)=0$ and so the $R$-module $\gf_S(M)$ is $S$-torsionfree. 
\item $\gf_S\gf_S(M) = \gf_S(M)$. 
\end{enumerate}
\end{lemma}

{\it Proof}. 1. $\gt_S\gf_S(M)=  \gt_S(\bM /\tor_{\bS}(\bM ))= \tor_{\bS}(\bM /\tor_{\bS}(\bM ) ) =0$.

2. $\gf_S\gf_S(M) \simeq \gf_S(M)/ \gt_S\gf_S(M)= \gf_S(M)$ since $\gt_S\gf_S(M)=0$, by statement 1.  $\Box $\\

Theorem \ref{20Jan19} is a criterion for the functor $S^{-1}:M\ra S^{-1}M$ to be exact. 

\begin{theorem}\label{20Jan19}
Let $R$ be a ring, $S\in \mLsRa$  where $*\in \{ l, \emptyset\}$, $\bR = R/ \ga $ and $\bS:= (S+\ga ) / \ga $. The functor $S^{-1}$ is exact iff for all $R$-modules $M_1$ and $M_2$ such that $M_1\subseteq M_2$, the $\bR$-modules $M_1\cap \ga M_2/\ga M_1$ is $\bS$-torsion. 
\end{theorem}

{\it Proof}. The theorem follows from the exact sequence (\ref{Slcm1}). 
 $\Box $

\begin{corollary}\label{b20Jan19}
Let $R$ be a ring, $S\in \mLsRa$ where $*\in \{ l, \emptyset\}$, $\bR = R/ \ga $ and $\bS:= (S+\ga ) / \ga $. If the functor $S^{-1}$ is exact then the $\bR$-module $\ga / \ga^2$ is $\bS$-torsion (recall that $\bS\in \Den_l(\bR , 0)$, by Theorem \ref{L11Jan19}.(1)).
\end{corollary}

{\it Proof}. Applying Theorem \ref{20Jan19} to the pair of $R$-modules $M_1=\ga \subseteq M_2= R$, we conclude that the $\bR$-module $(M_1\cap \ga M_2) / \ga M_1=\ga/ \ga^2$ is $\bS$-torsion. $\Box$\\

For an $R$-module $M$ and $S\in \mLsRa$, we have a descending chain of $R$-modules 
$$ \gt_S(M)\supseteq \gt_S^2(M)\supseteq \cdots \supseteq \gt_S^n(M)\supseteq \cdots $$
where $\gt_S^n(M) = \gt_S\gt_S\cdots \gt_S (M)$, $n$ times.





\small{

Department of Pure Mathematics

University of Sheffield

Hicks Building

Sheffield S3 7RH

UK

email: v.bavula@sheffield.ac.uk}

\end{document}